\documentclass[preprint, 12pt]{article}

\usepackage{amsfonts,amsmath,amssymb}
\usepackage{hyperref}
\usepackage{subcaption}
\usepackage[section]{placeins}
\usepackage[utf8]{inputenc}
\usepackage{tikz}
\usepackage{lineno}

\definecolor{darkblue}{rgb}{0,0.08,0.45}
\hypersetup{colorlinks,linkcolor=darkblue,citecolor=darkblue,filecolor=darkblue,urlcolor=darkblue}
\hypersetup{bookmarksopen,bookmarksnumbered}

\newcommand{\dy}{\ensuremath{\,\mathrm{d}\mathbf{y}}}
\newcommand{\normal}{\ensuremath{\hat{\mathbf{n}}}}
\newcommand{\utot}{\ensuremath{u_{\mathrm{tot}}}}
\newcommand{\usca}{\ensuremath{u_{\mathrm{sca}}}}
\newcommand{\uinc}{\ensuremath{u_{\mathrm{inc}}}}
\newcommand{\SLP}{\ensuremath{\mathcal{V}}}
\newcommand{\DLP}{\ensuremath{\mathcal{K}}}

\newcommand{\SL}{\ensuremath{V}}
\newcommand{\DL}{\ensuremath{K}}
\newcommand{\AD}{\ensuremath{T}}
\newcommand{\HS}{\ensuremath{D}}
\newcommand{\ID}{\ensuremath{I}}
\newcommand{\CA}{\ensuremath{A}}
\newcommand{\NtD}{\ensuremath{\Lambda_\mathrm{NtD}}}
\newcommand{\DtN}{\ensuremath{\Lambda_\mathrm{DtN}}}

\newcommand{\NtDmei}{\ensuremath{\Lambda_{\mathrm{NtD},m}^\pm}}
\newcommand{\DtNmei}{\ensuremath{\Lambda_{\mathrm{DtN},m}^\pm}}
\newcommand{\OsrcNtD}{\ensuremath{L_\mathrm{NtD}}}
\newcommand{\OsrcDtN}{\ensuremath{L_\mathrm{DtN}}}
\newcommand{\OsrcNtDm}{\ensuremath{L_{\mathrm{NtD},m}}}
\newcommand{\OsrcDtNm}{\ensuremath{L_{\mathrm{DtN},m}}}
\newcommand{\OsrcNtDme}{\ensuremath{L_{\mathrm{NtD},m}^+}}
\newcommand{\OsrcDtNme}{\ensuremath{L_{\mathrm{DtN},m}^+}}
\newcommand{\OsrcNtDmi}{\ensuremath{L_{\mathrm{NtD},m}^-}}
\newcommand{\OsrcDtNmi}{\ensuremath{L_{\mathrm{DtN},m}^-}}

\newcommand{\traceD}{\ensuremath{\gamma_D}}
\newcommand{\traceN}{\ensuremath{\gamma_N}}

\newcommand{\traceDme}{\ensuremath{\gamma_{D,m}^+}}
\newcommand{\traceDne}{\ensuremath{\gamma_{D,n}^+}}
\newcommand{\traceDmi}{\ensuremath{\gamma_{D,m}^-}}
\newcommand{\traceDmei}{\ensuremath{\gamma_{D,m}^\pm}}

\newcommand{\traceNme}{\ensuremath{\gamma_{N,m}^+}}
\newcommand{\traceNne}{\ensuremath{\gamma_{N,n}^+}}
\newcommand{\traceNmi}{\ensuremath{\gamma_{N,m}^-}}
\newcommand{\traceNmei}{\ensuremath{\gamma_{N,m}^\pm}}

\title{Frequency-robust preconditioning of boundary integral equations for acoustic transmission\footnote{© 2022. This manuscript version is made available under the CC-BY-NC-ND 4.0 license. This manuscript is published in the Journal of Computational Physics in final form at \url{https://doi.org/10.1016/j.jcp.2022.111229}.}}
\author{Elwin van 't Wout\thanks{Institute for Mathematical and Computational Engineering, School of Engineering and Faculty of Mathematics, Pontificia Universidad Católica de Chile, Santiago, Chile. Contact: e.wout@uc.cl} \and 
Seyyed R.~Haqshenas\thanks{Department of Mechanical Engineering, University College London, London, United Kingdom.} \and
Pierre Gélat\footnotemark[3] \and
Timo Betcke\thanks{Department of Mathematics, University College London, London, United Kingdom.} \and
Nader Saffari\footnotemark[3]}
\date{\today}

\begin{document}

\maketitle

\begin{abstract}
	The scattering and transmission of harmonic acoustic waves at a penetrable material are commonly modelled by a set of Helmholtz equations. This system of partial differential equations can be rewritten into boundary integral equations defined at the surface of the objects and solved with the boundary element method (BEM). High frequencies or geometrical details require a fine surface mesh, which increases the number of degrees of freedom in the weak formulation. Then, matrix compression techniques need to be combined with iterative linear solvers to limit the computational footprint. Moreover, the convergence of the iterative linear solvers often depends on the frequency of the wave field and the objects' characteristic size. Here, the robust PMCHWT formulation is used to solve the acoustic transmission problem. An operator preconditioner based on on-surface radiation conditions (OSRC) is designed that yields frequency-robust convergence characteristics. Computational benchmarks compare the performance of this novel preconditioned formulation with other preconditioners and boundary integral formulations. The OSRC preconditioned PMCHWT formulation effectively simulates large-scale problems of engineering interest, such as focused ultrasound treatment of osteoid osteoma.
\end{abstract}

\section{Introduction}

The simulation of scattering and transmission of harmonic acoustic waves at penetrable objects is of great interest to a wide range of engineering problems, such as biomedical ultrasound, underwater surveillance, and auralisation, among others~\cite{lahaye2017modern}. These physical phenomena, and many more, can be modelled with the Helmholtz equation in the volume and boundary conditions at the material interfaces. A wide range of numerical algorithms is available to approximate the solution of the boundary value problem directly, including the finite element and finite difference method~\cite{bergman2018computational}. Another approach is to consider integral formulations of the boundary value problem and solve it with volume integral equations~\cite{costabel2015spectrum} or boundary integral equations~\cite{marburg2018boundary}.

When scattering into an unbounded exterior domain is considered, the boundary element method (BEM) is often the preferred methodology~\cite{nedelec2001acoustic, steinbach2008numerical, sauter2010boundary}. Volumetric discretisation techniques, such as the finite element method, require artificial boundary conditions to limit the volumetric space under consideration. Techniques for the approximation of the outgoing waves include perfectly matched layers~\cite{berenger1994perfectly}, absorbing boundary conditions~\cite{engquist1977absorbing}, and infinite elements~\cite{gerdes1998summary}, among others. Differently, boundary integral formulations automatically satisfy the Sommerfeld radiation condition that models outgoing waves. Hence, no artificial boundary conditions are necessary, and fields can be evaluated at any point in space through Green's function and representation formulas for the surface potentials~\cite{marburg2008computational}.

For acoustic transmission into penetrable materials, transmission conditions couple exterior with interior fields and guarantee continuity of the pressure field and normal particle velocity. Boundary integral formulations can model transmission through different material regions when each is homogeneous. However, when a subdomain has heterogeneous material parameters, no Green's functions are available anymore. In that case, volumetric solvers should be used and coupled to the boundary integral formulations for the homogeneous subdomains~\cite{johnson1980coupling}. Here, we will consider homogeneous materials and use BEM-BEM coupling at the interfaces.

The BEM has several numerical properties that make it an efficient algorithm for acoustic wave propagation. Firstly, since Green's functions are explicitly used, the numerical algorithm is practically devoid of dispersion and dissipation errors, and a fixed number of elements per wavelength suffice for accurate results~\cite{graham2015wen, marburg2002six}. Secondly, only interfaces need to be meshed, simplifying the often cumbersome mesh generation procedure at large-scale geometries. Thirdly, the number of degrees of freedom does not depend on the white space between multiple scatterers~\cite{antoine2010computational}. These characteristics yield a quadratic scaling of the number of degrees of freedom with respect to the characteristic sizes of the scatterers.

A drawback of the BEM is that Green's functions lead to global operators and dense matrix arithmetic. In contrast, volumetric methods use sparse arithmetic. Therefore, the standard BEM has a computational complexity of $\mathcal{O}(n^2)$ with $n$ the number of degrees of freedom. Large-scale simulations are only feasible when accelerators are used, such as the fast-multipole method~\cite{greengard1987fast}, hierarchical matrix compression~\cite{borm2010efficient} or fast-Fourier transforms~\cite{bleszynski1996aim}. Most of these algorithms have a computational complexity of $\mathcal{O}(n \log(n))$. Furthermore, high-level open-source implementations of the BEM are available, such as BEMPP~\cite{smigaj2015solving}.

The driving frequency strongly influences the computational efficiency of the BEM. The number of degrees of freedom scales quadratically with the frequency since a fixed number of elements per wavelength is sufficient to obtain accurate results with the BEM. Furthermore, most acceleration algorithms such as fast-multipole and hierarchical matrix compression techniques have frequency-dependent performance~\cite{engquist2018approximate}. Moreover, the dense matrix equation is typically solved with iterative linear solvers~\cite{marburg2003performance}, whose convergence often deteriorates when increasing the frequency~\cite{galkowski2019wavenumber} or material contrast~\cite{wout2022highcontrast}. The actual convergence of the linear solvers depends on the choice of boundary integral formulation and the design of the preconditioner~\cite{antoine2021introduction}. One of the most widely used formulations for penetrable objects is the Poggio-Miller-Chang-Harrington-Wu-Tsai (PMCHWT) formulation~\cite{poggio1973integral}. For rigid scatterers, the OSRC preconditioner is very efficient at high frequencies~\cite{darbas2013combining, wout2015fast}. Here, we will design an OSRC preconditioner for the PMCHWT formulation for multiple penetrable objects.

Preconditioning of the dense set of linear equations is often necessary to perform large-scale simulations with the BEM. Among the many preconditioners, one can distinguish between algebraic and operator preconditioning. The algebraic preconditioners use properties of the discretisation matrix and include incomplete LU decompositions~\cite{sakuma2008fast} and the sparse approximate inverse~\cite{carpentieri2005combining}, among others. Differently, the operator preconditioners consider the characteristics of the continuous formulation~\cite{hiptmair2006operator, kirby2010functional}, such as mapping properties and Calderón projection formulas~\cite{steinbach1998construction}. Examples are mass, Calderón and OSRC preconditioners. Since operator preconditioners take the specific characteristics of boundary integral formulations into account, they often outperform algebraic preconditioners~\cite{antoine2008integral}. Furthermore, they do not require access to the discretisation matrix and are, therefore, feasible combinations with any matrix compression technique~\cite{darbas2013combining}. The OSRC preconditioner is one of the most effective techniques for large-scale simulations at high frequencies~\cite{wout2015fast}.

The OSRC operators are approximations of the Dirichlet-to-Neumann (DtN) map that links Dirichlet and Neumann traces of the acoustic field~\cite{antoine2008advances}. Since OSRC operators are local approximations of the DtN map, they are much quicker to solve than full boundary integral equations, but at the expense of inaccuracies in the solution~\cite{antoine2001fast}. Many different OSRC operators have been designed, based on wave decoupling techniques~\cite{moore1988theory}, absorbing boundary conditions~\cite{bayliss1980radiation}, and pseudo-differential operators~\cite{antoine1999bayliss}, among others. Furthermore, OSRC operators serve different purposes. Firstly, they are artificial boundary conditions that truncate exterior domains for the finite element and finite difference method~\cite{tsynkov1998numerical}. Secondly, they approximate scattering at objects~\cite{kriegsmann1987new}. Thirdly, they are efficient coupling operators in combined field integral equations such as the Burton-Miller and Brakhage-Werner formulations~\cite{antoine2004analytic}. The latter is a preconditioning technique since an OSRC coupling operator results in well-conditioned discretisation matrices of the BEM~\cite{antoine2005alternative, antoine2007generalized}. Finally, OSRC operators are also available for electromagnetics~\cite{darbas2006generalized, el2014approximate} and elastodynamics~\cite{chaillat2015approximate, darbas2015well}.

The OSRC operators approximate outgoing waves and, therefore, only model direct reflection and no multiple scattering~\cite{antoine2008advances}. This model limitation complicates the simulation of non-convex surfaces, multiple objects, and transmission. Nevertheless, in practice, OSRC preconditioning works well for objects with cavities~\cite{darbas2013combining}. For multiple scattering between objects, separate OSRC operators work well as preconditioner~\cite{wout2015fast}. However, when used to solve scattering phenomena directly, the OSRC approximations at each surface need to be coupled into a global system of equations~\cite{acosta2015surface, alzubaidi2016formulation}. In the case of penetrable domains, OSRC operators have been designed for generalized impedance boundary conditions for thin coatings~\cite{antoine2006improved}, as a weak coupling between FEM and BEM~\cite{caudron2020optimized} but not yet for BEM-BEM coupling.

In this study, we use OSRC operators as left preconditioners of the boundary integral formulation. This approach is different to the standard OSRC preconditioning, where the OSRC operator is a coupling parameter in a combined field integral equation. We apply OSRC preconditioning to transmission and multiple scattering and show its effectiveness in computational benchmarks. To the best of the authors' knowledge, the only application of OSRC preconditioning to the BEM for multiple scattering and transmission problems is in our previous work, where it was applied to the simulation of high-intensity focused ultrasound~\cite{haqshenas2021fast} and as part of a large benchmarking exercise~\cite{wout2021benchmarking}. Here, we will provide the mathematical foundation of this technique and convergence studies compared to other formulations.

The numerical benchmarks show that solving the BEM system with the GMRES algorithm requires more iterations when frequency increases. The OSRC preconditioner slows down this efficiency deterioration and is robust at high frequency ranges. Here, we define a preconditioner to be \emph{frequency-robust} when the speedup gain of preconditioning improves with frequency increase. The speedup gain is defined as the time to solve the preconditioned formulation divided by the solution time for the standard (mass-matrix preconditioned) PMCHWT formulation.

The formulation of the BEM and the OSRC preconditioner will be derived in Section~\ref{sec:formulations}. Computational benchmarks with different formulations and preconditioners will be presented in Section~\ref{sec:results}.

\section{Formulation}
\label{sec:formulations}

This study considers harmonic acoustic scattering and transmission at multiple penetrable homogeneous objects. This can be modelled by a set of Helmholtz equations in the homogeneous domains and boundary conditions at the interfaces of the materials. Then, boundary integral formulations are designed and solved with the BEM. The computational performance of the iterative linear method will be improved with preconditioning.

\subsection{Helmholtz equations}

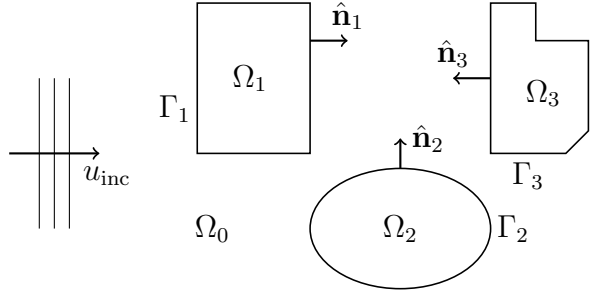
\begin{figure}[!ht]
	\centering
	\begin{tikzpicture}
		\draw[thin](-4.8,-1)--(-4.8,1);
		\draw[thin](-4.6,-1)--(-4.6,1);
		\draw[thin](-4.4,-1)--(-4.4,1);
		\draw[thick, ->](-5.2,0)--(-4,0);
		\node at (-3.9,-0.3) {$\uinc$};
		
		\node at (-2.5,-1) {$\Omega_0$};
		
		\draw[semithick] (-2.7,0) -- (-1.2,0) -- (-1.2,2) -- (-2.7,2) -- cycle;
		\node at (-2,1) {$\Omega_1$};
		\node at (-3,0.6) {$\Gamma_1$};
		\draw[thick, ->](-1.2,1.5)--(-0.7,1.5) node[above] {$\normal_1$};
		
		\draw[semithick](0,-1) ellipse (1.2 and 0.8) node {$\Omega_2$};
		\node at (1.5,-1) {$\Gamma_2$};
		\draw[thick, ->](0,-0.2)--(0,0.2) node[right] {$\normal_2$};
		
		\draw[semithick] (1.2,0) -- (2.2,0) -- (2.5,0.3) -- (2.5,1.5) -- (1.8,1.5) -- (1.8,2) -- (1.2,2) -- cycle;
		\node at (1.9,0.8) {$\Omega_3$};
		\node at (1.7,-0.3) {$\Gamma_3$};
		\draw[thick, ->](1.2,1)--(0.7,1) node[above] {$\normal_3$};
	\end{tikzpicture}
	\caption{A sketch of the geometry of the model.}
	\label{fig:domain}
\end{figure}

Let us consider multiple homogeneous domains denoted by $\Omega_m$ with $m=1,2,\dots,\ell$. These three-dimensional objects are all bounded and embedded in an unbounded exterior domain denoted by $\Omega_0$, as depicted in Figure~\ref{fig:domain}. The boundaries of the interior domains are denoted by $\Gamma_m$ for $m=1,2,\dots,\ell$ and $\Gamma_0$ is the union of all interfaces. The interfaces are assumed to be piecewise smooth such that unique normal vectors can be defined, which are denoted by $\normal_m$, have unit length, are directed towards the exterior domain, and correspond to interface $\Gamma_m$ for $m=1,2,\dots,\ell$, respectively.

The incident wave is given by a plane wave field with frequency~$f$, unit amplitude and direction~$\hat{\mathbf{d}}$. In the absence of scatterers, the acoustic field would be given by the incident wave field denoted by $\uinc$. The presence of objects generates a scattered field denoted by $\usca$ and the total field is given by $\utot = \uinc + \usca$. For the design of the boundary integral formulations, it is convenient to decompose the total field into separate fields associated to each subdomain. For this purpose, let us consider
\begin{equation}
	u_0(\mathbf{x}) =
	\begin{cases}
		\usca(\mathbf{x}) & \text{for } \mathbf{x} \in \Omega_0, \\
		-\uinc(\mathbf{x}) & \text{for } \mathbf{x} \notin \Omega_0
	\end{cases}
\end{equation}
the field associated to the exterior domain and
\begin{equation}
	u_m(\mathbf{x}) =
	\begin{cases}
		\utot(\mathbf{x}) & \text{for } \mathbf{x} \in \Omega_m, \\
		0 & \text{for } \mathbf{x} \notin \Omega_m
	\end{cases}
\end{equation}
for $m=1,2,\dots,\ell$ the fields associated to the different objects.

The equations of motion for the acoustic field are given by
\begin{equation}
	\label{eq:helmholtz}
	\begin{cases}
		-\Delta u_j - k_j u_j = 0 & \text{in } \Omega_j \text{ for } j=0,1,2,\dots,\ell; \\
		\traceDme (u_0 + \uinc) = \traceDmi u_m & \text{at } \Gamma_m \text{ for } m=1,2,\dots,\ell; \\
		\frac1{\rho_0} \traceNme (u_0 + \uinc) = \frac1{\rho_m} \traceNmi u_m & \text{at } \Gamma_m \text{ for } m=1,2,\dots,\ell; \\
		\lim_{\mathbf{r} \to \infty} |\mathbf{r}| (\partial_{|\mathbf{r}|} u_0 - ik_0 u_0) = 0.
	\end{cases}
\end{equation}
The first equation corresponds to the Helmholtz equation inside each subdomain where $k_j$ denotes the wavenumber. The second and third equations are interface conditions that model the continuity of the pressure field and normal particle velocity, respectively, where $\rho_j$ denotes the density of the material in~$\Omega_j$. The fourth equation is the Sommerfeld radiation condition that states that the scattered field radiates outward at infinity. Here, the trace operators are defined as
\begin{linenomath}
\begin{align}
	\traceDme f(\mathbf{x}) &= \lim_{\mathbf{y} \to \mathbf{x}} f(\mathbf{y}) && \text{for } \mathbf{x} \in \Gamma_m \text{ and } \mathbf{y} \in \Omega_0, \\
	\traceDmi f(\mathbf{x}) &= \lim_{\mathbf{y} \to \mathbf{x}} f(\mathbf{y}) && \text{for } \mathbf{x} \in \Gamma_m \text{ and } \mathbf{y} \in \Omega_m, \\
	\traceNme f(\mathbf{x}) &= \lim_{\mathbf{y} \to \mathbf{x}} \nabla_{\mathbf{y}} f(\mathbf{y}) \cdot \normal_m(\mathbf{x}) && \text{for } \mathbf{x} \in \Gamma_m \text{ and } \mathbf{y} \in \Omega_0, \\
	\traceNmi f(\mathbf{x}) &= \lim_{\mathbf{y} \to \mathbf{x}} \nabla_{\mathbf{y}} f(\mathbf{y}) \cdot \normal_m(\mathbf{x}) && \text{for } \mathbf{x} \in \Gamma_m \text{ and } \mathbf{y} \in \Omega_m
\end{align}
\end{linenomath}
for $m=1,2,\dots,\ell$ where the subindices $D$ and $N$ denote Dirichlet and Neumann traces, respectively; the subindex $m$ denotes the respective subdomain; and the superindices $+$ and $-$ denote exterior and interior traces, respectively. The operator $\nabla_{\mathbf{y}}$ denotes the gradient with respect to the variable $\mathbf{y}$. Notice that the normals are always pointing towards the exterior, regardless of the type of trace.

\subsection{Boundary integral formulations}

The Helmholtz system~\eqref{eq:helmholtz} is defined as a volumetric problem. Inside each subdomain, Green's functions are available:
\begin{equation}
	G_m(\mathbf{x},\mathbf{y}) = \frac{e^{ik_m |\mathbf{x} - \mathbf{y}|}}{4\pi |\mathbf{x} - \mathbf{y}|} \qquad \text{for } \mathbf{x},\mathbf{y} \in \Omega_m \text{ and } \mathbf{x} \ne \mathbf{y}
\end{equation}
for $m=0,1,2,\dots,\ell$ where $i$ denotes the imaginary unit. Hence, the field inside each subdomain can be defined by representation formulas~\cite{steinbach2008numerical} as
\begin{linenomath}
\begin{align}
	u_0 &= \sum_{n=1}^\ell \left(\SLP_{0,n} \tilde\psi_n^+ - \DLP_{0,n} \tilde\phi_n^+\right),
	\label{eq:representation:exterior:multiple} \\
	u_m &= \SLP_m \tilde\psi_m^- - \DLP_m \tilde\phi_m^- && \text{for } m=1,2,\dots,\ell;
	\label{eq:representation:interior:multiple}
\end{align}
\end{linenomath}
where the single-layer and double-layer potential integral operators are given by
\begin{linenomath}
\begin{align}
	[\SLP_{0,m} \psi](\mathbf{x}) &= \iint_{\Gamma_m} G_0(\mathbf{x},\mathbf{y}) \psi(\mathbf{y}) \dy && \text{for } \mathbf{x} \notin \Gamma_m, \\
	[\DLP_{0,m} \phi](\mathbf{x}) &= \iint_{\Gamma_m} \normal_m \cdot \nabla_{\mathbf{y}} \, G_0(\mathbf{x},\mathbf{y}) \phi(\mathbf{y}) \dy && \text{for } \mathbf{x} \notin \Gamma_m, \\
	[\SLP_m \psi](\mathbf{x}) &= \iint_{\Gamma_m} G_m(\mathbf{x},\mathbf{y}) \psi(\mathbf{y}) \dy && \text{for } \mathbf{x} \notin \Gamma_m, \\
	[\DLP_m \phi](\mathbf{x}) &= \iint_{\Gamma_m} \normal_m \cdot \nabla_{\mathbf{y}} \, G_m(\mathbf{x},\mathbf{y}) \phi(\mathbf{y}) \dy && \text{for } \mathbf{x} \notin \Gamma_m,
\end{align}
\end{linenomath}
for $m=1,2,\dots,\ell$.
Taking traces of the representation formulas yields
\begin{linenomath}
\begin{align}
	\traceDme u_0 &= \sum_{n=1}^\ell \SL_{0,mn} \tilde\psi_n^+ - \frac12 \tilde\phi_m^+ - \sum_{n=1}^\ell \DL_{0,mn} \tilde\phi_n^-, \\
	\traceNme u_0 &= -\frac12 \tilde\psi_m^+ + \sum_{n=1}^\ell \AD_{0,mn} \tilde\psi_n^+ + \sum_{n=1}^\ell \HS_{0,mn} \tilde\phi_n^+, \\
	\traceDmi u_m &= \SL_m \tilde\psi_m^- + \frac12 \tilde\phi_m^- - \DL_m \tilde\phi_m^-, \\
	\traceNmi u_m &= \frac12 \tilde\psi_m^- + \AD_m \tilde\psi_m^- + \HS_m \tilde\phi_m^-,
\end{align}
\end{linenomath}
for $m=1,2,\dots,\ell$, where 
\begin{linenomath}
\begin{align}
	\SL_m: H^{-\frac12}(\Gamma_m) &\to H^{\frac12}(\Gamma_m) \nonumber \\
	\psi &\mapsto \iint_{\Gamma_m} G_m (\cdot,\mathbf{y}) \psi(\mathbf{y}) \dy, \\
	\DL_m: H^{\frac12}(\Gamma_m) &\to H^{\frac12}(\Gamma_m) \nonumber \\
	\phi &\mapsto \iint_{\Gamma_m} \normal_m \cdot \nabla_{\mathbf{y}} \, G_m(\cdot,\mathbf{y}) \phi(\mathbf{y}) \dy, \\
	\AD_m: H^{-\frac12}(\Gamma_m) &\to H^{-\frac12}(\Gamma_m) \nonumber \\
	\psi &\mapsto \normal_m \cdot \nabla \iint_{\Gamma_m} G_m(\cdot,\mathbf{y}) \psi(\mathbf{y}) \dy, \\
	\HS_m: H^{\frac12}(\Gamma_m) &\to H^{-\frac12}(\Gamma_m) \nonumber \\
	\phi &\mapsto -\normal_m \cdot \nabla \iint_{\Gamma_m} \normal_m \cdot \nabla_{\mathbf{y}} \, G_m(\cdot,\mathbf{y}) \phi(\mathbf{y}) \dy,
\end{align}
\end{linenomath}
denote the single-layer, double-layer, adjoint double-layer, and hypersingular boundary integral operator, respectively; and
\begin{linenomath}
\begin{align}
	\SL_{0,mn}: H^{-\frac12}(\Gamma_n) &\to H^{\frac12}(\Gamma_m) \nonumber \\
	\psi &\mapsto \iint_{\Gamma_n} G_0 (\cdot,\mathbf{y}) \psi(\mathbf{y}) \dy, \\
	\DL_{0,mn}: H^{\frac12}(\Gamma_n) &\to H^{\frac12}(\Gamma_m) \nonumber \\
	\phi &\mapsto \iint_{\Gamma_n} \normal_n \cdot \nabla_{\mathbf{y}} \, G_0(\cdot,\mathbf{y}) \phi(\mathbf{y}) \dy, \\
	\AD_{0,mn}: H^{-\frac12}(\Gamma_n) &\to H^{-\frac12}(\Gamma_m) \nonumber \\
	\psi &\mapsto \normal_m \cdot \nabla \iint_{\Gamma_n} G_0(\cdot,\mathbf{y}) \psi(\mathbf{y}) \dy, \\
	\HS_{0,mn}: H^{\frac12}(\Gamma_n) &\to H^{-\frac12}(\Gamma_m) \nonumber \\
	\phi &\mapsto -\normal_m \cdot \nabla \iint_{\Gamma_n} \normal_n \cdot \nabla_{\mathbf{y}} \, G_0(\cdot,\mathbf{y}) \phi(\mathbf{y}) \dy,
\end{align}
\end{linenomath}
the boundary integral operators for the cross scattering in the exterior.
The surface potentials are traces of the total field~\cite{steinbach2008numerical}, specifically,
\begin{linenomath}
\begin{align}
	\tilde\phi_m^+ &= -\traceDme \utot, \\
	\tilde\phi_m^- &=  \traceDmi \utot, \\
	\tilde\psi_m^+ &= -\traceNme \utot, \\
	\tilde\psi_m^- &=  \traceNmi \utot
\end{align}
\end{linenomath}
for $m=1,2,\dots,\ell$.
Hence, the boundary integral equations satisfy
\begin{linenomath}
\begin{align}
	\left(\frac12\bar\ID_m + A_{0,mm} \right)
	\begin{bmatrix} \traceDme \utot \\ \traceNme \utot \end{bmatrix} + \sum_{n=1, n \ne m}^\ell A_{0,mn}
	\begin{bmatrix} \traceDne \utot \\ \traceNne \utot \end{bmatrix}
	&= \begin{bmatrix} \traceDme \uinc \\ \traceNme \uinc \end{bmatrix},
	\label{eq:calderon:exterior:traces} \\
	\left(-\frac12\bar\ID_m + A_m \right)
	\begin{bmatrix} \traceDmi \utot \\ \traceNmi \utot \end{bmatrix}
	&= \begin{bmatrix} 0 \\ 0 \end{bmatrix}
	\label{eq:calderon:interior:traces}
\end{align}
\end{linenomath}
for $m=1,2,\dots,\ell$ where $\bar\ID_m = \begin{bmatrix} \ID_m & 0 \\ 0 & \ID_m \end{bmatrix}$ with $\ID_m$ the identity operator at interface~$m$, and
\begin{equation}
	\CA_m = \begin{bmatrix} -\DL_m & \SL_m \\ \HS_m & \AD_m \end{bmatrix}
\end{equation}
denotes the Calderón boundary integral operator.

Now, let us define a unique pair of traces at each interface as
\begin{linenomath}
\begin{align}
	\phi_m &= \traceDme \utot = \traceDmi \utot, \\
	\psi_m &= \traceNme \utot = \frac{\rho_0}{\rho_m} \traceNmi \utot.
\end{align}
\end{linenomath}
Then,
\begin{linenomath}
\begin{align}
	\left(\frac12\bar\ID_m + A_{0,mm} \right)
	\begin{bmatrix} \phi_m \\ \psi_m \end{bmatrix} + \sum_{n=1, n \ne m}^\ell A_{0,mn}
	\begin{bmatrix} \phi_n \\ \psi_n \end{bmatrix}
	&= \begin{bmatrix} \traceDme \uinc \\ \traceNme \uinc \end{bmatrix}
	\label{eq:calderon:exterior:potentials} \\
	\left(-\frac12\bar\ID_m + \widehat{A}_m \right)
	\begin{bmatrix} \phi_m \\ \psi_m \end{bmatrix}
	&= \begin{bmatrix} 0 \\ 0 \end{bmatrix}
	\label{eq:calderon:interior:potentials}
\end{align}
\end{linenomath}
for $m=1,2,\dots,\ell$ where
\begin{equation}
	\widehat{\CA}_m = \begin{bmatrix} -\DL_m & \frac{\rho_m}{\rho_0} \SL_m \\ \frac{\rho_0}{\rho_m} \HS_m & \AD_m \end{bmatrix}.
\end{equation}
Finally, unique boundary integral formulations can be obtained by taking linear combinations of the boundary integral equations~\cite{mitzner1966acoustic}. For example, taking the sum of Eqns.~\eqref{eq:calderon:exterior:potentials} and~\eqref{eq:calderon:interior:potentials} yields
\begin{equation}
	\label{eq:pmchwt}
	\left(A_{0,mm} + \widehat{A}_m\right)
	\begin{bmatrix} \phi_m \\ \psi_m \end{bmatrix} + \sum_{n=1, n \ne m}^\ell A_{0,mn}
	\begin{bmatrix} \phi_n \\ \psi_n \end{bmatrix}
	= \begin{bmatrix} \traceDme \uinc \\ \traceNme \uinc \end{bmatrix}
\end{equation}
for $m=1,2,\dots,\ell$ which is the PMCHWT formulation~\cite{poggio1973integral, chang1974surface, wu1977scattering-bor}.
Alternatively, taking the difference of Eqns.~\eqref{eq:calderon:exterior:potentials} and~\eqref{eq:calderon:interior:potentials} yields
\begin{equation}
	\label{eq:muller}
	\left(\ID_m + A_{0,mm} - \widehat{A}_m\right)
	\begin{bmatrix} \phi_m \\ \psi_m \end{bmatrix} + \sum_{n=1, n \ne m}^\ell A_{0,mn}
	\begin{bmatrix} \phi_n \\ \psi_n \end{bmatrix}
	= \begin{bmatrix} \traceDme \uinc \\ \traceNme \uinc \end{bmatrix}
\end{equation}
for $m=1,2,\dots,\ell$ which is the Müller formulation~\cite{muller1957grundprobleme}.
The PMCHWT formulation is a first-kind boundary integral equation and the Müller formulation is of second-kind.

\subsection{Boundary element method}

The weak form of the boundary integral equations are discretised with a Galerkin method~\cite{sauter2010boundary}. At the material interfaces, a triangular surface mesh will be used. The discrete function space for $H^{\frac12}(\Gamma_m)$ will be piecewise linear functions (P1) associated to the vertices in the mesh. For $H^{-\frac12}(\Gamma_m)$, one can use either piecewise constant functions (P0) associated to the triangles in the mesh, piecewise constant functions on the dual mesh (that is, node-based patches on the barycentric refined mesh~\cite{buffa2007dual}), or P1 functions. Here, P1 functions will always be used.

\subsection{Mass preconditioning}

Let us write the discrete set of equations as
\begin{equation}
	\label{eq:weakform}
	B \mathbf{x} = \mathbf{f}
\end{equation}
where $B$ denotes the discretisation matrix, $\mathbf{x}$ the unknown coefficients of the basis functions, and $\mathbf{f}$ the traces of the incident wave field.
Alternatively, mass preconditioning yields the strong form
\begin{equation}
	\label{eq:strongform}
	M^{-1} B \mathbf{x} = M^{-1} \mathbf{f}
\end{equation}
that maps the domain space to the range space when $M$ is the identity operator mapping the range space to the dual-to-range space~\cite{betcke2020product}. Notice that when a combination of P0 and P1 elements is used, this mapping is not well-defined. In that case, either dual functions or P1 functions need to be used.

\subsection{Calderón preconditioning}

The Calderón operator is a projection~\cite{steinbach2008numerical}, that is, $A_m^2 = \ID_m$ for $m=0,1,2,\dots,\ell$. This property can be used to create operator preconditioners for the PMCHWT formulation in the form of
\begin{equation}
	M_2^{-1} C M_1^{-1} B \mathbf{x} = M_2^{-1} C M_1^{-1} \mathbf{f}.
\end{equation}
Here, we will consider the block-diagonal Calderón preconditioner
\begin{equation}
	\label{eq:calderon}
	C = \begin{bmatrix} A_{0,11} + \widehat{A}_1 && \emptyset \\ & \ddots \\ \emptyset && A_{0,\ell\ell} + \widehat{A}_\ell \end{bmatrix}.
\end{equation}
Notice that other versions can be designed as well~\cite{yan2010comparative, cools2011calderon, wout2021benchmarking}. This preconditioner does not require any additional assembly time since the operators are already available in the model. In each GMRES iteration, the preconditioner step involves the matrix-vector multiplication of $C$, which has $2\ell$ Calderón operators.

\subsection{OSRC preconditioning}

The Neumann-to-Dirichlet (NtD) and Dirichlet-to-Neumann (DtN) maps, also known as Steklov-Poincaré and Poincaré-Steklov operators~\cite{sauter2010boundary}, are implicitly defined as
\begin{linenomath}
\begin{align}
	\NtDmei: H^{-\frac12}(\Gamma_m) &\to H^{\frac12}(\Gamma_m) \nonumber \\
	\traceNmei\utot &\mapsto \traceDmei\utot, \\
	\DtNmei: H^{\frac12}(\Gamma_m) &\to H^{-\frac12}(\Gamma_m) \nonumber \\
	\traceDmei\utot &\mapsto \traceNmei\utot,
\end{align}
\end{linenomath}
where the superscript $\pm$ denotes the exterior and interior versions, respectively. Notice that these maps involve inverse boundary integral operators. Hence, no closed-form expressions are available for general surfaces and numerical approximations involve the solution of a system of boundary integral operators. This solution procedure is too expensive for the direct application of NtD and DtN maps in boundary integral formulations. Instead, accurate approximations of these maps can be obtained with on-surface radiation conditions. Among the many different definitions of OSRC operators, the ones based on pseudo-differential equations are among the most accurate for general convex objects~\cite{antoine2008advances, alzubaidi2016formulation}, which are defined as
\begin{linenomath}
\begin{align}
	\OsrcNtDme &= \frac1{ik_0} \left(\ID_m + \frac{\Delta_{\Gamma_m}}{k_{\epsilon,0,m}^2}\right)^{-\frac12}, \\
	\OsrcNtDmi &= \frac1{ik_m} \left(\ID_m + \frac{\Delta_{\Gamma_m}}{k_{\epsilon,m}^2}\right)^{-\frac12}, \\
	\OsrcDtNme &= ik_0 \left(\ID_m + \frac{\Delta_{\Gamma_m}}{k_{\epsilon,0,m}^2}\right)^{\frac12}, \\
	\OsrcDtNmi &= ik_m \left(\ID_m + \frac{\Delta_{\Gamma_m}}{k_{\epsilon,m}^2}\right)^{\frac12}
\end{align}
\end{linenomath}
where $\Delta_{\Gamma_m}$ denotes the Laplace-Beltrami operator on surface~$\Gamma_m$. Furthermore,
\begin{linenomath}
\begin{align*}
	\OsrcNtD^\pm&: H^{-\frac12}(\Gamma) \to H^{\frac12}(\Gamma),
	\\
	\OsrcDtN^\pm&: H^{\frac12}(\Gamma) \to H^{-\frac12}(\Gamma).
\end{align*}
\end{linenomath}
Singularities are avoided by using a damped wavenumber $k_{\epsilon,0,m} \ne k_0$ and $k_{\epsilon,m} \ne k_m$ for $m=1,2,\dots,\ell$. This artificial damping actually improves its accuracy at high frequencies~\cite{antoine2005improved, antoine2006improved}.
In practice, a rule of thumb for choosing the damping is $k_{\epsilon,0,m} = k_0 (1 + 0.4 i (k_0 R_m)^{-\frac23})$ and $k_{\epsilon,m} = k_m (1 + 0.4 i (k_m R_m)^{-\frac23})$ where $R_m$ denotes the radius of the object~$\Omega_m$. This choice is optimal for a sphere~\cite{antoine2005alternative} and extensions exist that use the local curvature of the surface~\cite{antoine2008advances}. For simplicity, this study uses the default expression.

In the case of a single scatterer, one has
\begin{linenomath}
\begin{align}
	\NtD^\pm: H^{\frac12}(\Gamma) &\to H^{\frac12}(\Gamma), \nonumber \\
	\HS^\pm \traceD^\pm u^\pm &\mapsto \left(\mp\tfrac12\ID - \DL^\pm\right) \traceD^\pm u^\pm, \\
	\DtN^\pm: H^{-\frac12}(\Gamma) &\to H^{-\frac12}(\Gamma), \nonumber \\
	\SL^\pm \traceN^\pm u^\pm &\mapsto \left(\mp\tfrac12\ID + \AD^\pm\right) \traceN^\pm u^\pm
\end{align}
\end{linenomath}
for any interior solution~$u^-$ and radiating solution~$u^+$ of the Helmholtz equation~\cite{steinbach2008numerical}. In the case of the Helmholtz transmission system~\eqref{eq:helmholtz}, the interior Calderón system~\eqref{eq:calderon:interior:traces} yields an equivalent result, that is,
\begin{linenomath}
\begin{align}
	(\NtD^- \HS_m) \traceDmi \utot &= \left(\tfrac12\ID_m - \DL_m\right) \traceDmi \utot, \\
	(\DtN^- \SL_m) \traceNmi \utot &= \left(\tfrac12\ID_m + \AD_m\right) \traceNmi \utot
\end{align}
\end{linenomath}
for $m=1,2,\dots,\ell$. This shows that using the NtD and DtN operators are good preconditioners of the hypersingular and single-layer operators, respectively, since the right-hand side consists of second-kind boundary integral operators that are well-conditioned. In practice, the OSRC operators are used instead of the NtD and DtN operators. Still, the operator product will be well-conditioned since this change only adds a compact perturbation to the right-hand side~\cite{darbas2013combining}. For the exterior domain, the situation is more complicated due to the presence of the incident wave field as well as the multiple reflections of the scattered field. Exterior OSRC operators for multiple scattering can be designed~\cite{acosta2015surface, alzubaidi2016formulation} but require explicit cross-interactions between objects and, therefore, require more computational overhead. For this reason, we will consider independent OSRC operators at each surface. Even though multiple scattering is not included in these operators, computational experiments show an excellent performance for preconditioning at rigid scatterers~\cite{wout2015fast} and the benchmarks in this study confirm its effectiveness for transmission problems as well.

\subsection{The OSRC preconditioned PMCHWT formulation}

The PMCHWT formulation~\eqref{eq:pmchwt} has single-layer and hypersingular operators on the off-diagonal blocks of the discretisation matrix. Hence, a permuted version will be used and combined with the OSRC preconditioner as follows
\begin{linenomath}
\begin{align}
	&\begin{bmatrix} \OsrcDtNm & 0 \\ 0 & \OsrcNtDm \end{bmatrix}
	\begin{bmatrix} \SL_{0,11} + \frac{\rho_1}{\rho_0} \SL_1 & -\DL_{0,11} - \DL_1 \\ \AD_{0,11} + \AD_1 & \HS_{0,11} + \frac{\rho_0}{\rho_1} \HS_1 \end{bmatrix}
	\begin{bmatrix} \psi_1 \\ \phi_1 \end{bmatrix} \nonumber \\
	&\qquad = \begin{bmatrix} \OsrcDtNm & 0 \\ 0 & \OsrcNtDm \end{bmatrix}
	\begin{bmatrix} \traceDme \uinc \\ \traceNme \uinc \end{bmatrix}
\end{align}
\end{linenomath}
in the case of a single scatterer. The preconditioning can readily be extended to multiple scattering as
\begin{linenomath}
\begin{align}
	&\begin{bmatrix} L_1 & 0 & \cdots & 0 \\ 0 & L_2 & \cdots & 0 \\ \vdots & \vdots & \ddots & \vdots \\ 0 & 0 & \cdots & L_m \end{bmatrix}
	\begin{bmatrix}
		A_{0,11}^\star + \widehat{A}_1^\star & A_{0,12}^\star & \cdots & A_{0,1\ell}^\star \\
		A_{0,21}^\star & A_{0,22}^\star + \widehat{A}_2^\star & \cdots & A_{0,2\ell}^\star \\
		\vdots & \vdots & \ddots & \vdots \\
		A_{0,\ell1}^\star & A_{0,\ell2}^\star & \cdots & A_{0,\ell\ell}^\star + \widehat{A}_\ell^\star \\
	\end{bmatrix}
	\begin{bmatrix} \varphi_1 \\ \varphi_2 \\ \vdots \\ \varphi_\ell \end{bmatrix} \nonumber \\
	&\qquad = \begin{bmatrix} L_1 & 0 & \cdots & 0 \\ 0 & L_2 & \cdots & 0 \\ \vdots & \vdots & \ddots & \vdots \\ 0 & 0 & \cdots & L_\ell \end{bmatrix}
	\begin{bmatrix} f_1 \\ f_2 \\ \vdots \\ f_\ell \end{bmatrix}
	\label{eq:pmchwt:osrc}
\end{align}
\end{linenomath}
with
\begin{linenomath}
\begin{align*}
	L_m &= \begin{bmatrix} \OsrcDtNm & 0 \\ 0 & \OsrcNtDm \end{bmatrix}, &
	\varphi_m &= \begin{bmatrix} \psi_m \\ \phi_m \end{bmatrix}, &
	f_m &= \begin{bmatrix} \traceDme \uinc \\ \traceNme \uinc \end{bmatrix},
\end{align*}
\end{linenomath}
and where the stars denote permutation, that is,
\begin{equation}
	\CA_m^\star = \begin{bmatrix} \SL_m & -\DL_m \\ \AD_m & \HS_m \end{bmatrix}
\end{equation}
the permuted Calderón operator. Notice that one could also use a block OSRC preconditioner given by
\begin{equation}
	L_m^\star = \begin{bmatrix} 0 & \OsrcNtDm \\ \OsrcDtNm & 0 \end{bmatrix}
\end{equation}
for the original (non-permuted) PMCHWT formulation.

The preconditioned formulation~\eqref{eq:pmchwt:osrc} uses a single OSRC operator for the sum of the exterior and interior single-layer and hypersingular operators. Hence, either the interior or the exterior wavenumber can be used for the OSRC preconditioner. In fact, any wavenumber could be used and the choice can be different at each interface. Taking a different wavenumber for the OSRC operator than for the Calderón operator is not an issue since a change in wavenumber results in a compact perturbation of the boundary integral operators~\cite{antoine2008integral, claeys2013multi, boubendir2015integral}.

The square-root in the OSRC operator needs to be localised for efficiency purposes~\cite{moore1988theory}. Here, a truncated Padé series expansion of $N_\text{Padé}$ terms and a branch cut with an angle~$\theta$ is used~\cite{antoine2008advances}. This results in a coupled set of $N_\text{Padé}$ equations each involving a surface Helmholtz equation with different complex-valued wavenumbers. The weak formulation of these equations can readily be solved by the Galerkin method with P1 elements as test and basis functions. Hence, at each surface two sets of $N_\text{Padé}$ linear equations need to be solved at each iteration of the linear solver. Since all these systems are sparse, a sparse LU decomposition is calculated once at the start of the linear solver and used in each iteration.

\section{Results}
\label{sec:results}

The novel application of OSRC preconditioning to the PMCHWT formulation promises to be a frequency-robust solver for transmission problems. In this section, the results of computational benchmarks will be shown, as well as the performance of the preconditioned formulation at large-scale geometries.

\subsection{Benchmark framework}

The benchmarks will consider the following boundary integral formulations:
\begin{enumerate}
	\item The PMCHWT formulation~\eqref{eq:pmchwt}.
	\item The Müller formulation~\eqref{eq:muller}.
	\item The PMCHWT formulation with the block-diagonal Calderón preconditioner~\eqref{eq:calderon}.
	\item The PMCHWT formulation with the OSRC preconditioner~\eqref{eq:pmchwt:osrc}.
\end{enumerate}
The strong form of all operators will be considered, in other words, mass preconditioning~\eqref{eq:strongform} is applied by default. The discretisation of the formulations is given by a Galerkin method with P1 elements for all test and basis functions.

For the OSRC operators, the parameters are chosen, if not stated otherwise, as $N_\text{Padé} = 4$ and $\theta = \pi/3$. The parameter $R_m$ in the damped wavenumber normally denotes the radius of the object. However, computational experience suggests that taking $R_m$ to be one tenth of the radius improves the effectiveness of OSRC preconditioning for transmission problems. The wavenumber for the OSRC preconditioner is either the interior or exterior wavenumber at each surface.

Version 4.3 of the library Gmsh~\cite{geuzaine2009gmsh} was used to generate the triangular surface meshes at the interfaces, with a mesh width of $h$. That is, each triangular element has a characteristic length smaller than $h$. The mesh width at interface~$m$ was chosen as
\begin{equation}
	h_m = \frac{\min\{ \lambda_0 , \lambda_m \}}{n_h}
\end{equation}
for $m=1,2,\dots,\ell$ where $\lambda_m$ is the wavelength in $\Omega_m$. In other words, there are at least $n_h$ elements per wavelength at both sides of each interface.

The discrete operators will be assembled with hierarchical matrix compression, with a tolerance of $10^{-5}$. As iterative linear solver, the GMRES algorithm is used, with a tolerance of $10^{-7}$ for the preconditioned residual and without restart.

The boundary integral formulations have been implemented with version~3 of the BEMPP library~\cite{smigaj2015solving, scroggs2017software}. The GMRES solver comes from the library SciPy version 1.2.1~\cite{scipy}. All simulations have been performed with hyperthreading activated, on a workstation with 16 processor cores (Intel\textregistered~Xeon(R) CPU E5-2683 v4\textcircled{a}2.10 GHz) and 512 GB RAM of shared memory.

Different material types will be used in the benchmarking with a frequency power law model for attenuation~\cite{hamilton1998nonlinear}. That is,
\begin{equation}
	k_m =  \frac{2\pi f}{c_m} + i \alpha_m (f \cdot 10^{-6})^{b_m}
\end{equation}
where $\alpha$ denotes the attenuation coefficient in Np~m$^{-1}$~Hz$^{-1}$ and $b$ an exponent. See Table~\ref{table:parameters:physical} for characteristic values of materials commonly found in acoustical engineering.

\begin{table}[!ht]
	\caption{Physical parameters of the scattering media~\cite{duck1990physical,itis2018}.}
	\label{table:parameters:physical}
	\centering
	\begin{tabular}{lrrrr}
		\hline\hline
		material & $\rho$ & $c$ & $\alpha$ & $b$ \\
		\hline
		water & 1000 & 1500 & 0.015 & 2 \\
		fat & 917 & 1412 & 9.334 & 1 \\
		bone & 1912 & 4080 & 47.20 & 1 \\
		\hline\hline
	\end{tabular}
\end{table}

\subsection{Model accuracy}

\begin{figure}[!ht]
	\centering
	\includegraphics[width=\columnwidth]{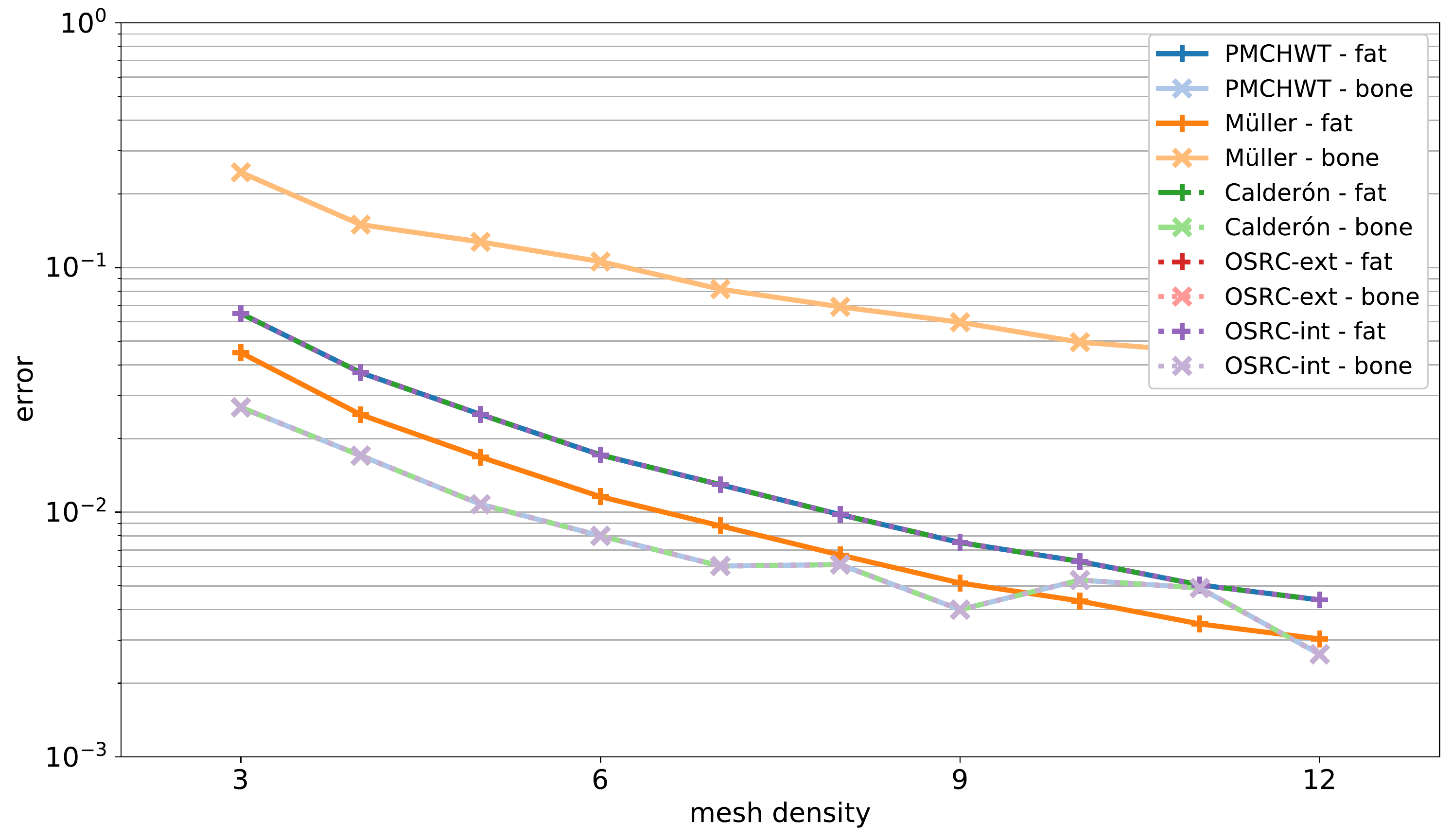}
	\caption{The error of the BEM relative to the analytical solution at a penetrable sphere, with respect to the mesh density. The error is measured in the $\ell_2$ norm of the amplitude of the pressure field evaluated on a grid of $101 \times 101$ points uniformly located on a square of size $3 \times 3$~cm in the $x$-$y$ plane and centered in the global origin. The incident wave field is a plane wave with frequency of 1~MHz.}
	\label{fig:error}
\end{figure}

Let us first test the accuracy of the formulations at a single spherical object, for which analytical solutions are available in the form of a series expansion in spherical harmonics. The pressure field obtained from the BEM has been calculated in the interior and exterior of the sphere, where the exterior has acoustic parameters resembling water and the interior models either fat or bone material. Figure~\ref{fig:error} presents the error of the field relative to the analytical solution. The preconditioned versions of the PMCHWT did not have noticeable differences in error measure compared to the PMCHWT formulation itself. The accuracy of the fields improves exponentially with the mesh density and six elements per wavelength is sufficient to obtain an accuracy of 10\% for the Müller formulation and 2\% for the preconditioned PMCHWT formulations.

\FloatBarrier
\subsection{Frequency dependency}

The frequency of the wave field has a profound influence on the computational performance of the BEM. The OSRC preconditioner is a sparse operator based on high-frequency approximation and is, therefore, expected to work well at high frequencies. To validate the computational efficiency, let us consider the same geometry as before: a sphere with radius 5~mm and a mesh with six elements per wavelength. This leads to a number of 1805, 7088, 15\,823, 27\,237, and 42\,619 nodes in the triangular surface mesh for 500~kHZ, 1~MHz, 1.5~MHz, 2~MHz, and 2.5~MHz, respectively.

\begin{figure}[!ht]
	\centering
	\begin{subfigure}[b]{\columnwidth}
		\centering
		\includegraphics[width=\columnwidth]{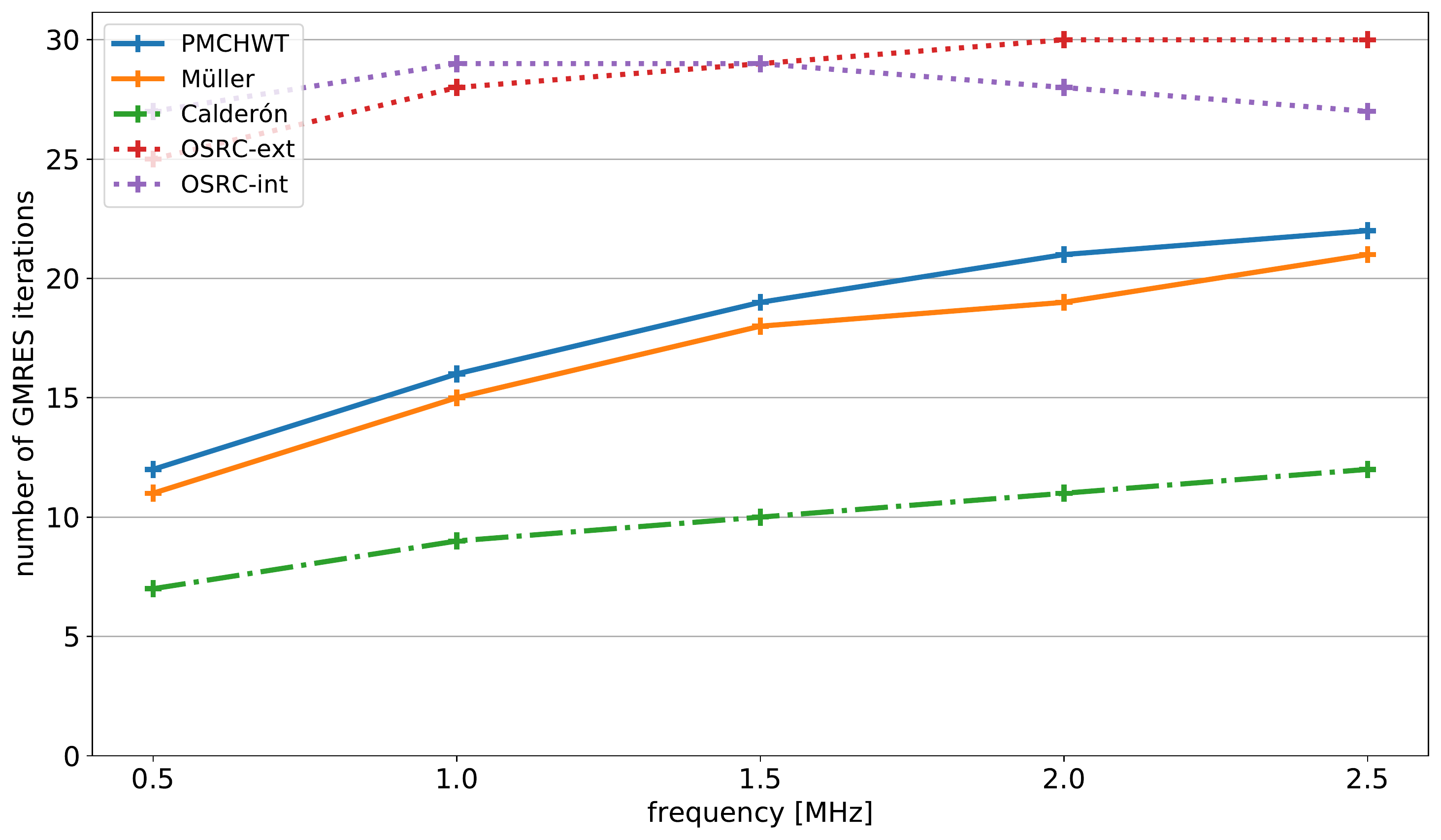}
		\caption{Exterior water, interior fat material.}
	\end{subfigure}
	\begin{subfigure}[b]{\columnwidth}
		\centering
		\includegraphics[width=\columnwidth]{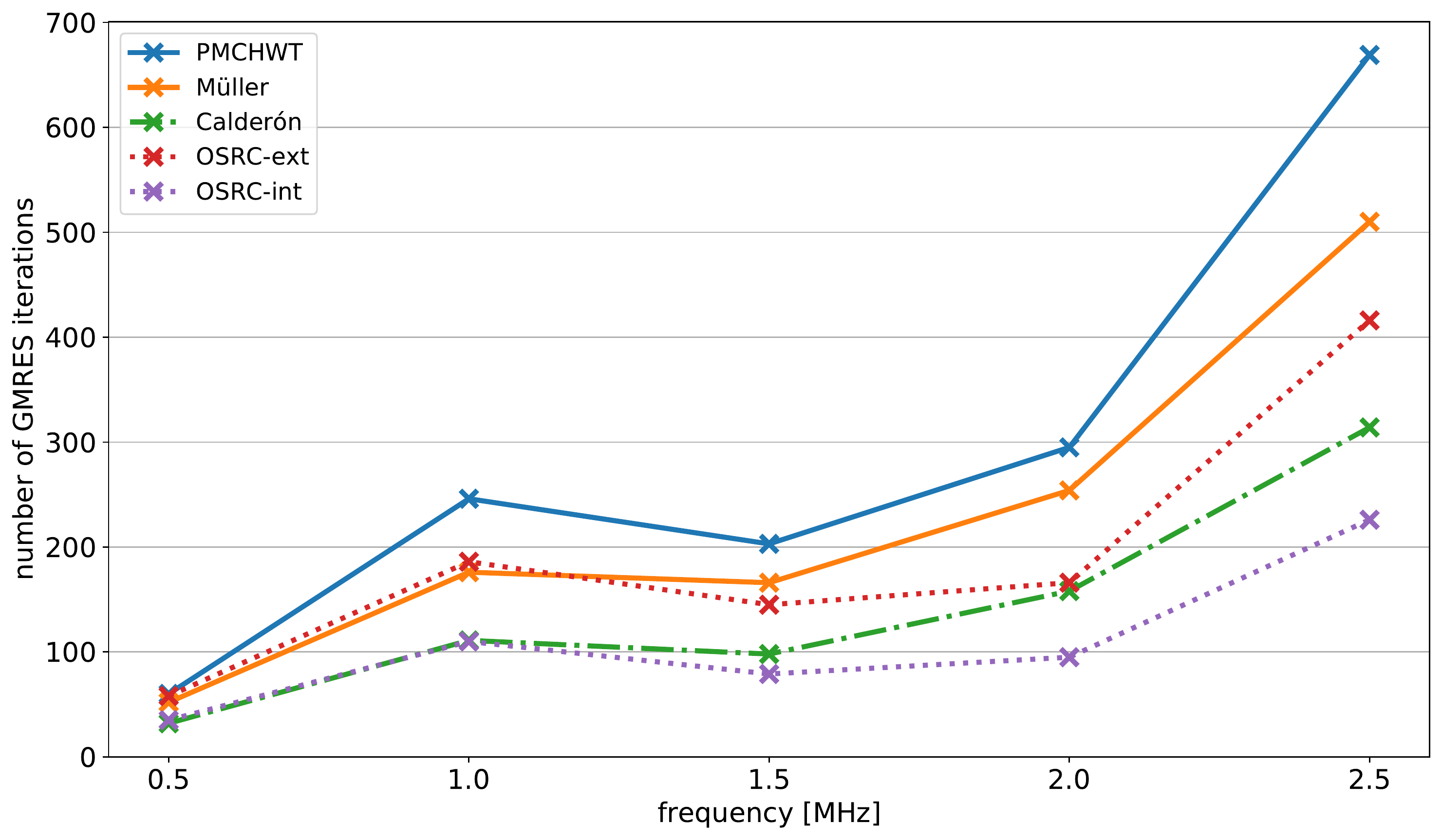}
		\caption{Exterior water, interior bone material.}
	\end{subfigure}
	\caption{The number of GMRES iterations with respect to the frequency. The geometry is a penetrable sphere with radius 5~mm.}
	\label{fig:frequency:iterations}
\end{figure}

Figure~\ref{fig:frequency:iterations} presents the convergence of the GMRES algorithm for increasing frequency. In the case of a water-fat interface, there is a low contrast in density and wavespeed between the materials. The PMCHWT and Müller formulations already converge very quickly for this situation and OSRC preconditioning does not reduce the number of GMRES iterations. Remember that the OSRC preconditioner is applied to the permuted PMCHWT formulation, not to the standard PMCHWT formulation. For this reason, the OSRC-PMCHWT formulation can be worse conditioned than the PMCHWT formulation. The tendency is a deterioration of efficiency for the (Calderón preconditioned) PMCHWT and Müller formulations while the number of GMRES iterations for the OSRC preconditioned PMCHWT remains constant with frequency. The situation is different for the water-bone interface, where higher material contrasts are present and the system is worse conditioned. The PMCHWT and Müller formulations require a large number of GMRES iterations, which can be reduced considerably with Calderón and OSRC preconditioning. The OSRC preconditioner with a wavenumber taken from the interior material performs best.

\begin{figure}[!ht]
	\centering
	\begin{subfigure}[b]{\columnwidth}
		\centering
		\includegraphics[width=\columnwidth]{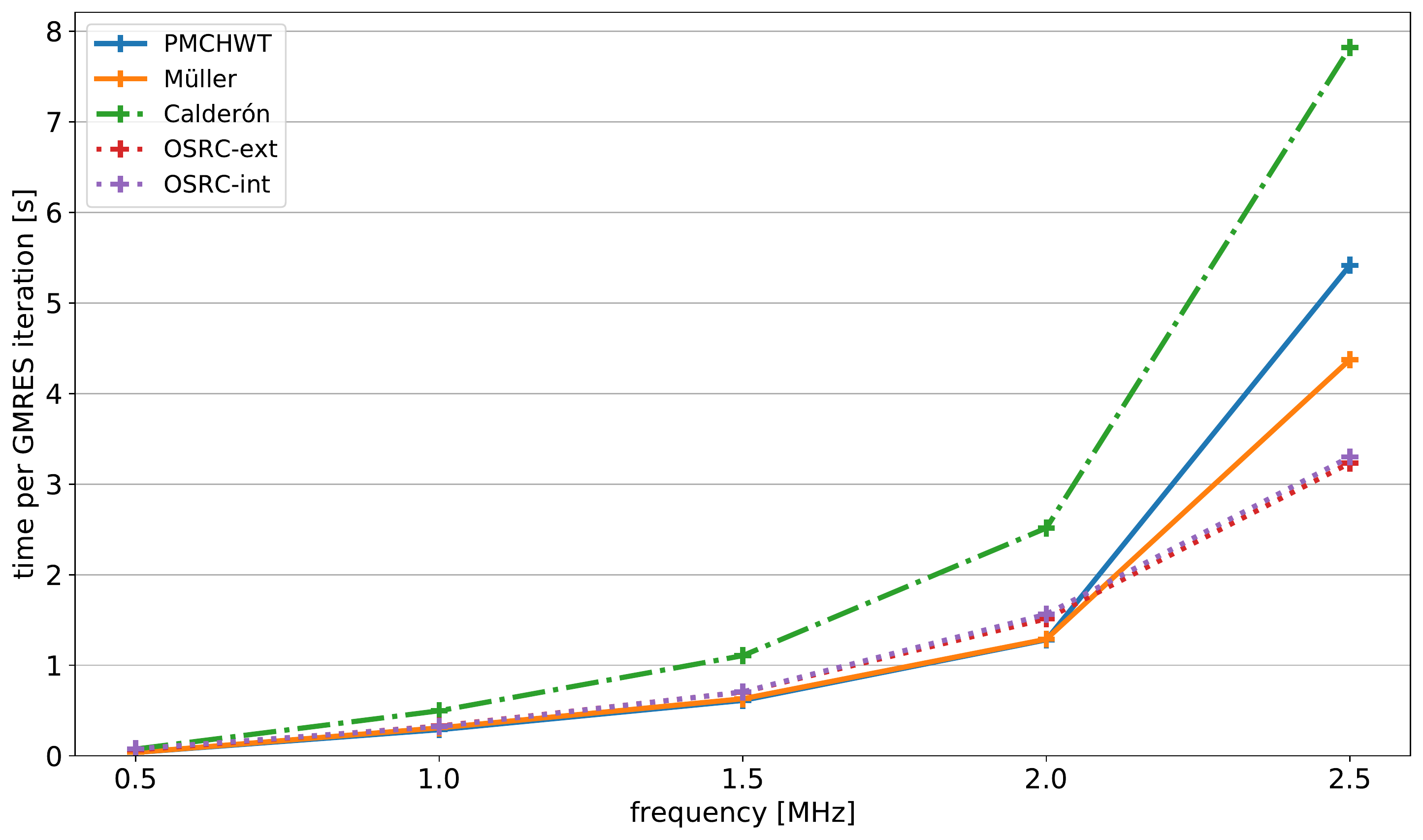}
		\caption{Exterior water, interior fat material.}
	\end{subfigure}
	\begin{subfigure}[b]{\columnwidth}
		\centering
		\includegraphics[width=\columnwidth]{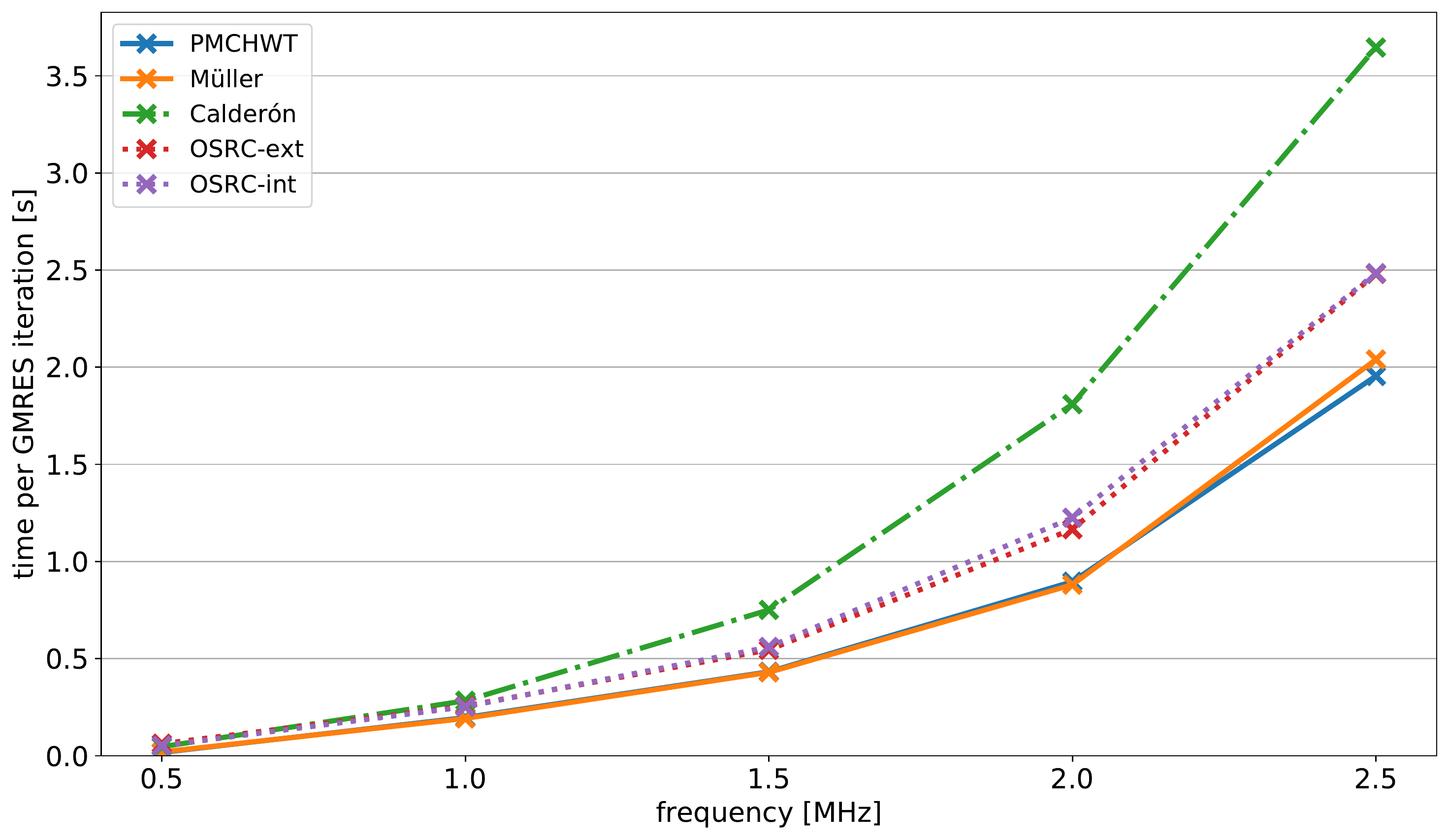}
		\caption{Exterior water, interior bone material.}
	\end{subfigure}
	\caption{The computation time (wall-clock time) per GMRES iteration with respect to the frequency. The geometry is a penetrable sphere with radius 5~mm.}
	\label{fig:frequency:time:iteration}
\end{figure}

\begin{figure}[!ht]
	\centering
	\begin{subfigure}[b]{\columnwidth}
		\centering
		\includegraphics[width=\columnwidth]{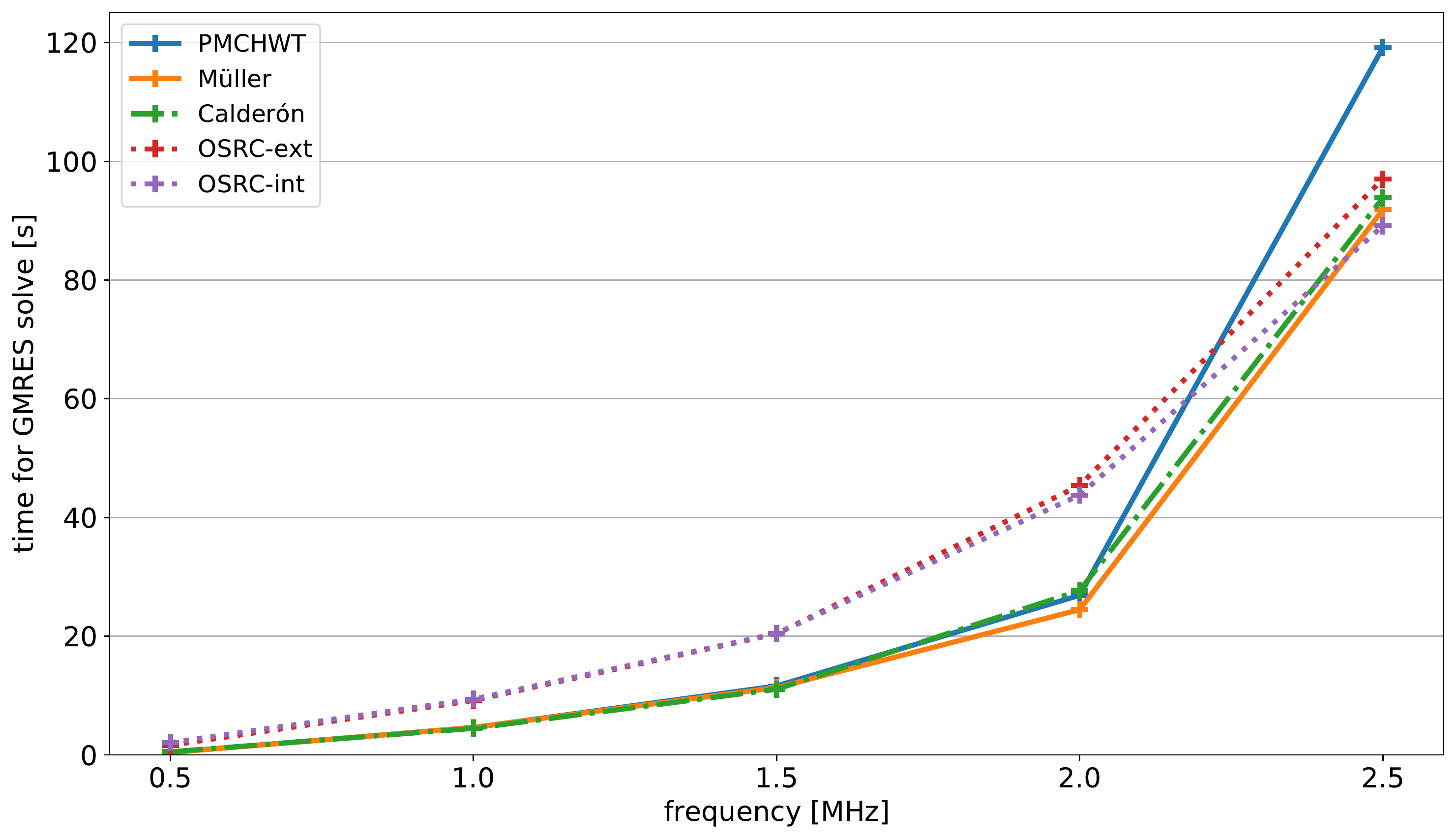}
		\caption{Exterior water, interior fat material.}
	\end{subfigure}
	\begin{subfigure}[b]{\columnwidth}
		\centering
		\includegraphics[width=\columnwidth]{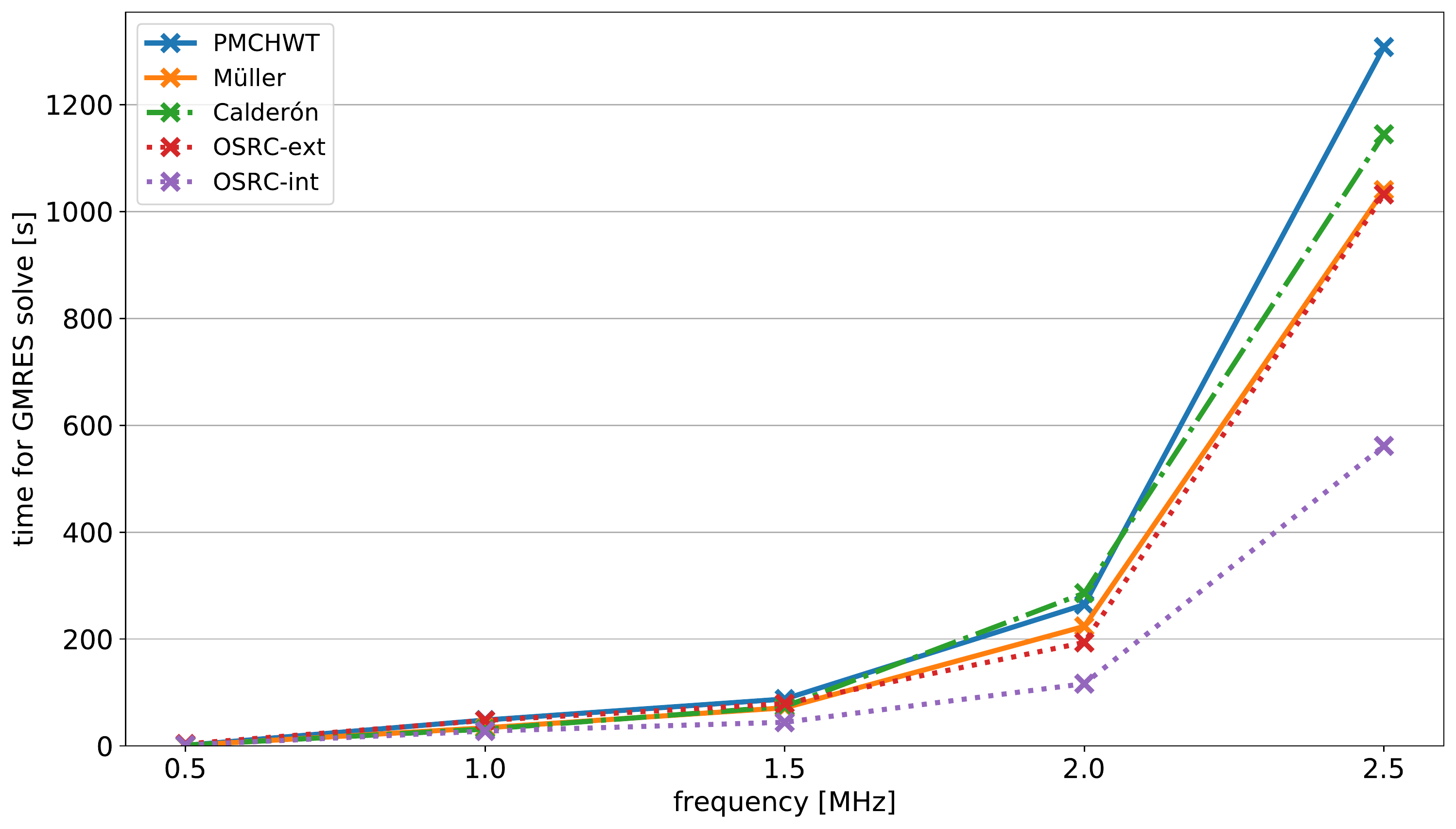}
		\caption{Exterior water, interior bone material.}
	\end{subfigure}
	\caption{The total time (wall-clock time) for the GMRES algorithm with respect to the frequency. The geometry is a penetrable sphere with radius 5~mm.}
	\label{fig:frequency:time:gmres}
\end{figure}

The computation time to assemble the matrices is the same for all formulations, since all of them require the interior and exterior Calderón operators. The computational overhead of creating the OSRC operators is always less than 0.5\% of the time to build the model matrix. Concerning the time for the GMRES algorithm, the most expensive part in each iteration is the matrix-vector multiplication of the system matrix, stored in compressed format. The other expensive operation for GMRES is the preconditioner step. Calderón preconditioning doubles the timing since the square of the matrix is being used. OSRC preconditioning requires the solution of a set of $N_\text{Padé}$ (chosen to be four) surface Helmholtz system, for which the sparse LU factorisation was calculated during the matrix assembly. This has little overhead and the OSRC preconditioned PMCHWT formulation is the best performing formulation for the challenging configurations with high frequencies and material contrasts, as can be seen in Figures~\ref{fig:frequency:time:iteration} and~\ref{fig:frequency:time:gmres}.

\FloatBarrier
\subsection{Domains with corners}

\begin{figure}[!ht]
	\centering
	\includegraphics[width=.8\columnwidth]{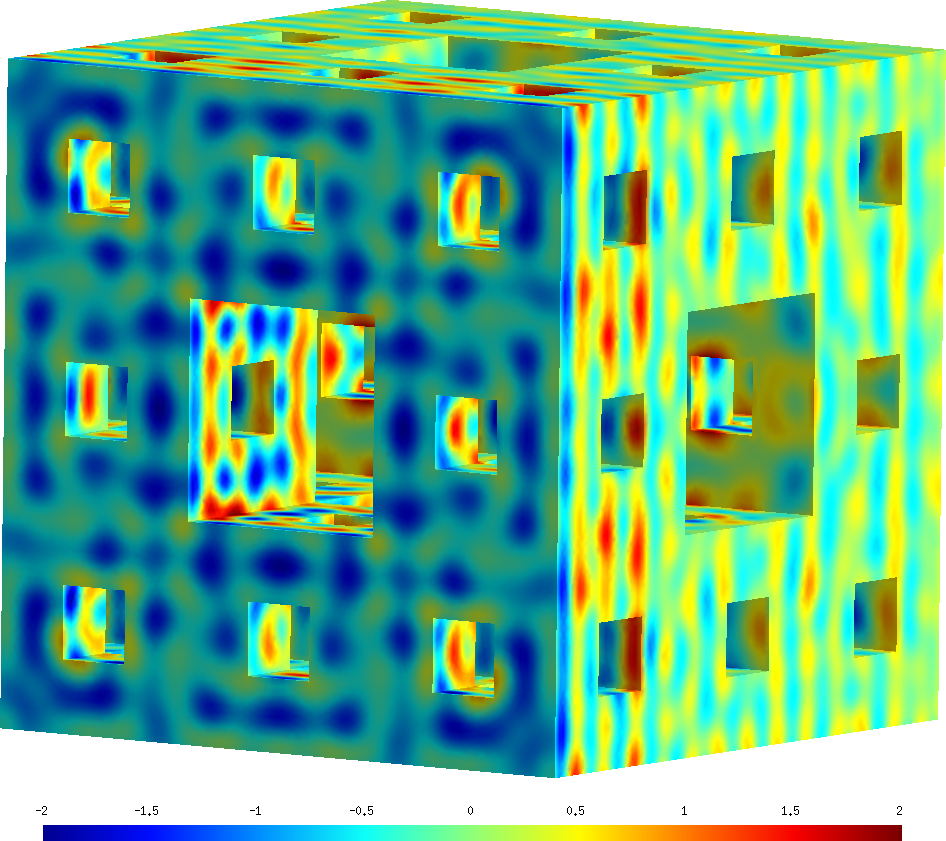}
	\caption{The acoustic pressure at the surface of a Menger sponge of level~2 with edges of 10~mm and the materials resemble water and bone for the exterior and interior domains, respectively. The incident wave field is a plane wave travelling towards the left face of the cube, with unit amplitude, and a frequency of 2~MHz.}
	\label{fig:sponge:geometry}
\end{figure}

\begin{figure}[!ht]
	\centering
	\begin{subfigure}[b]{\columnwidth}
		\centering
		\includegraphics[width=\columnwidth]{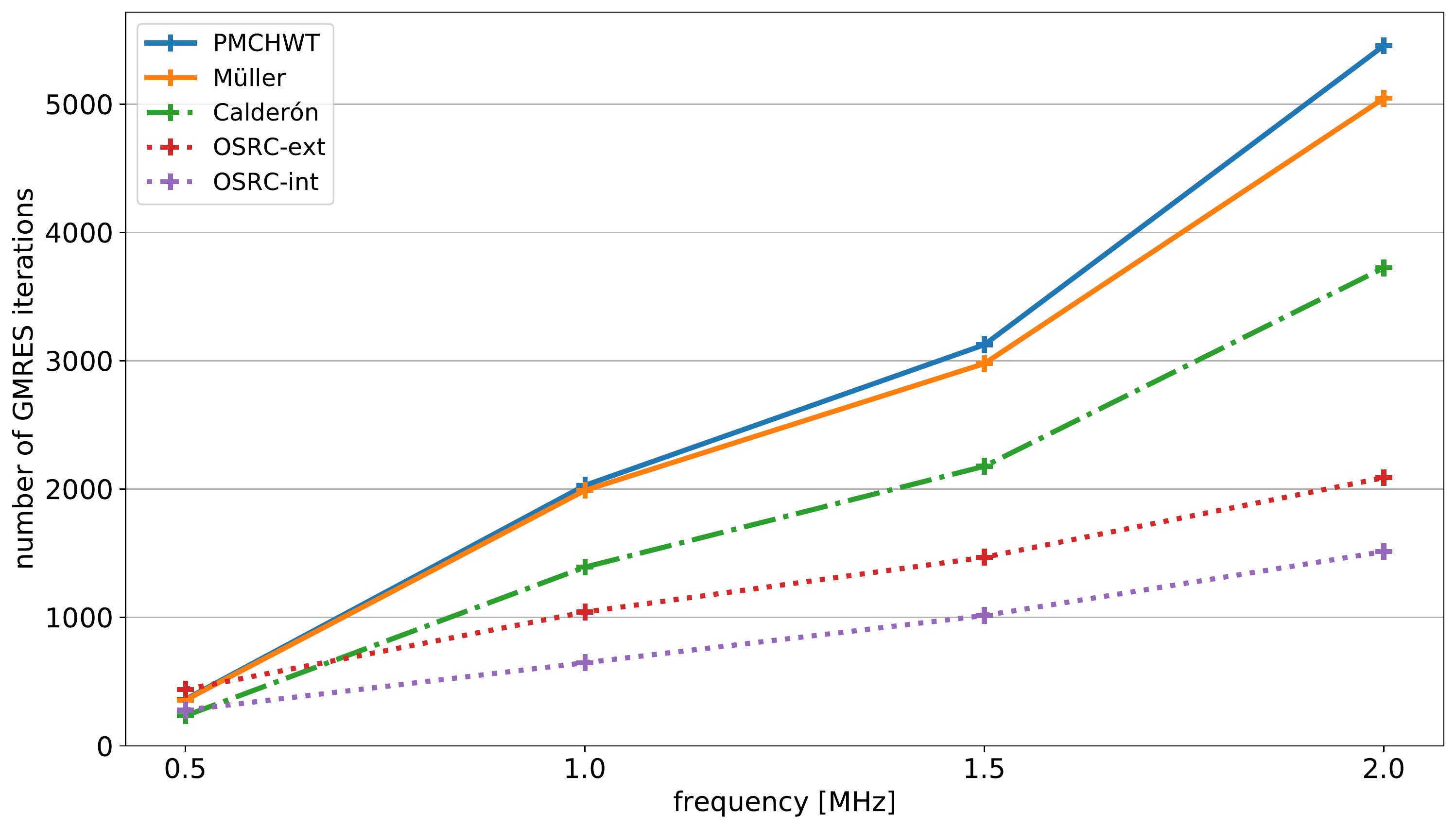}
		\caption{The number of GMRES iterations.}
	\end{subfigure}
	\begin{subfigure}[b]{\columnwidth}
		\centering
		\includegraphics[width=\columnwidth]{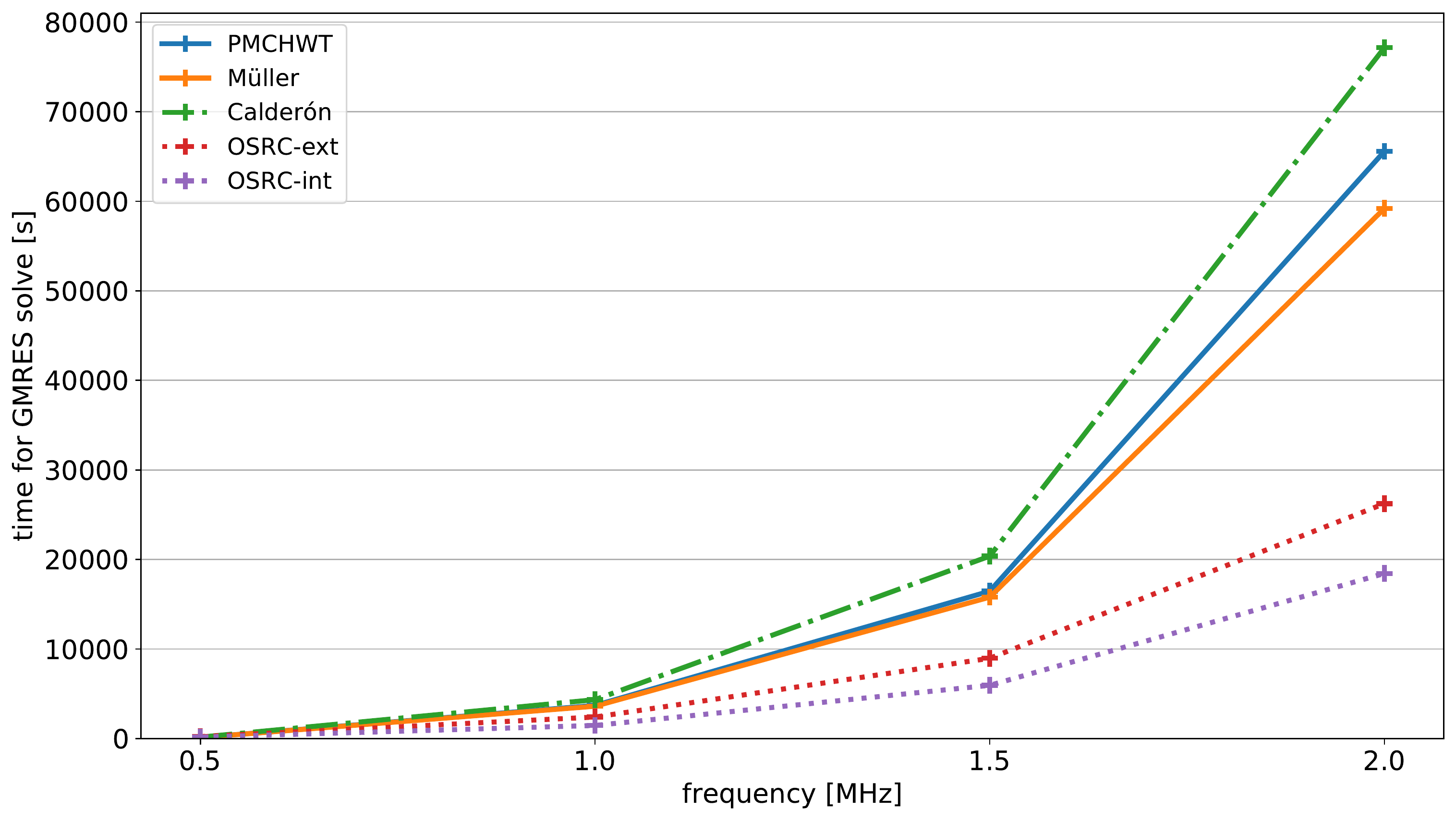}
		\caption{The total time (wall-clock time) for the GMRES algorithm.}
	\end{subfigure}
	\caption{The computational performance with respect to the frequency. The geometry is a Menger sponge of level~2 with edges of 10~mm and the materials resemble water and bone for the exterior and interior domains, respectively.}
	\label{fig:sponge:frequency}
\end{figure}

The OSRC operators are accurate at smooth domains~\cite{darbas2013combining} and the benchmark on a sphere confirmed the efficiency of OSRC preconditioners for the PMCHWT formulation. To test the feasibility of OSRC preconditioning on nonsmooth domains, let us consider a Menger sponge, which is a fractal volume with increasingly small inclusions. The outer dimension of the cube is 10\,mm for each edge and two levels of fractal divisions are considered, see Fig.~\ref{fig:sponge:geometry}. The GMRES convergence and compute time are measured at the frequencies 500~kHz, 1~MHz, 1.5~MHz, and 2~MHz with at least six mesh elements per wavelength, resulting in surface grids with 13\,568, 34\,975, 62\,291, and 103\,987~nodes, respectively. The computational results are presented in Fig.~\ref{fig:sponge:frequency}. The OSRC preconditioned PMCHWT formulation significantly outperforms the Müller, PMCHWT, and Calderón preconditioned PMCHWT in terms of GMRES convergence and solution time. For example, the OSRC preconditioner with interior wavenumber is more than three times faster than any of the standard formulations at 2~MHz. Furthermore, the speedup gained with OSRC preconditioning of the PMCHWT formulation improves with frequency. Comparing the 1~MHz with 2~MHz simulations, the speedup factor increases from 1.5 to 2.5 for exterior OSRC preconditioning and from 2.5 to 3.6 for interior OSRC, respectively. Hence, the OSRC preconditioner scales favorably with respect to frequency, compared to the other formulations. This computational benchmark corroborates the robustness of the OSRC preconditioner for challenging simulations that involve high-frequency wave fields and geometries with corners and inclusions.

\FloatBarrier
\subsection{Multiple scattering}

All boundary integral formulations and preconditioners considered in this study can be applied to multiple scattering at disjoint penetrable objects. Let us consider a configuration of four spheres with radius 5~mm with centers at locations ($\pm$15, $\pm$7.5, 0)~mm. Two of which resemble fat material and the other two bone. The surface mesh has six elements per wavelength, yielding a number of 6823, 26\,801, 59\,758, and 102\,549 nodes for a frequency of 500~kHz, 1~MHz, 1.5~MHz, and 2~MHz, respectively.

\begin{figure}[!ht]
	\centering
	\begin{subfigure}[b]{\columnwidth}
		\centering
		\includegraphics[height=.20\paperheight]{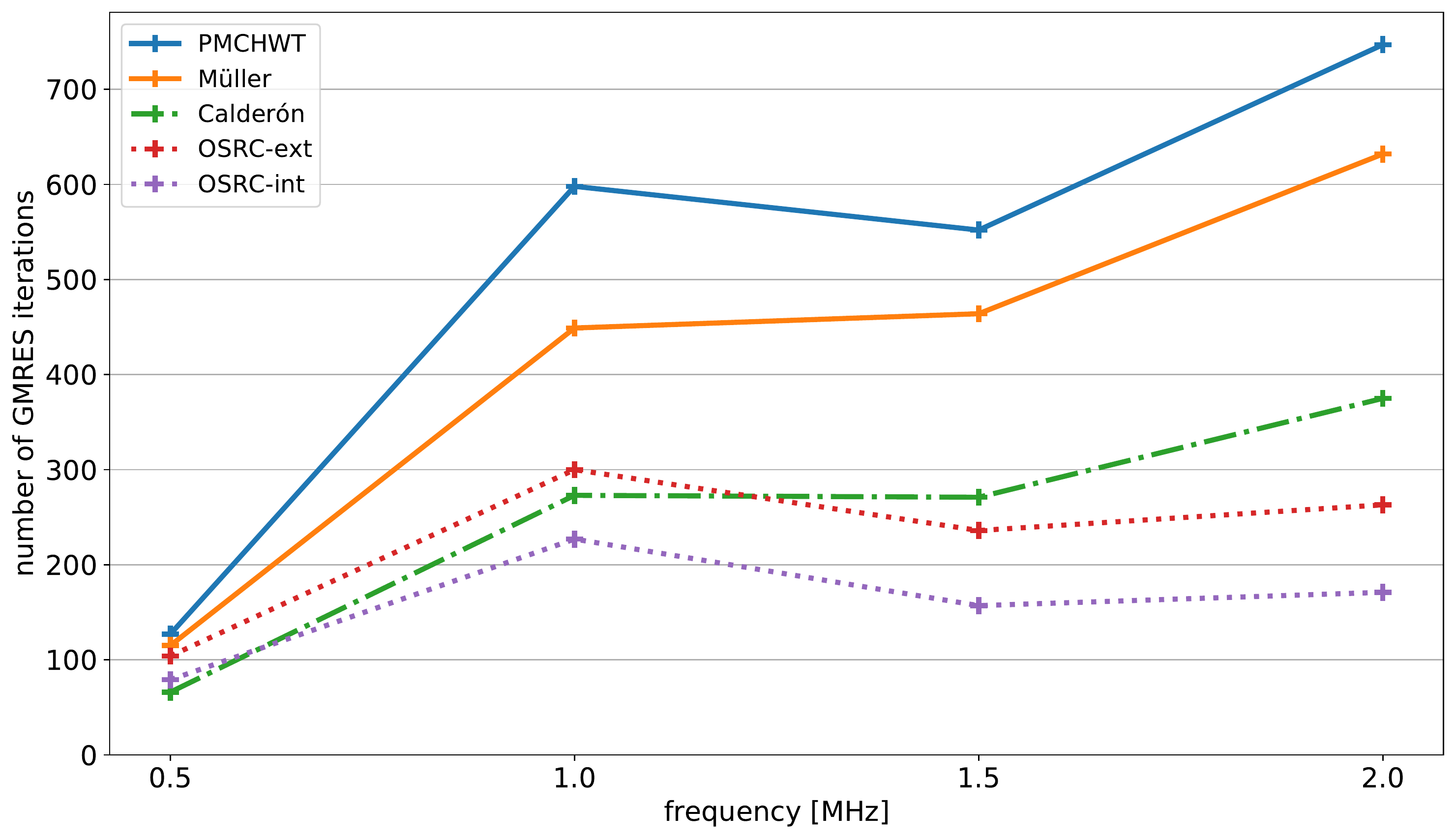}
		\caption{The number of GMRES iterations.}
	\end{subfigure}
	\begin{subfigure}[b]{\columnwidth}
		\centering
		\includegraphics[height=.20\paperheight]{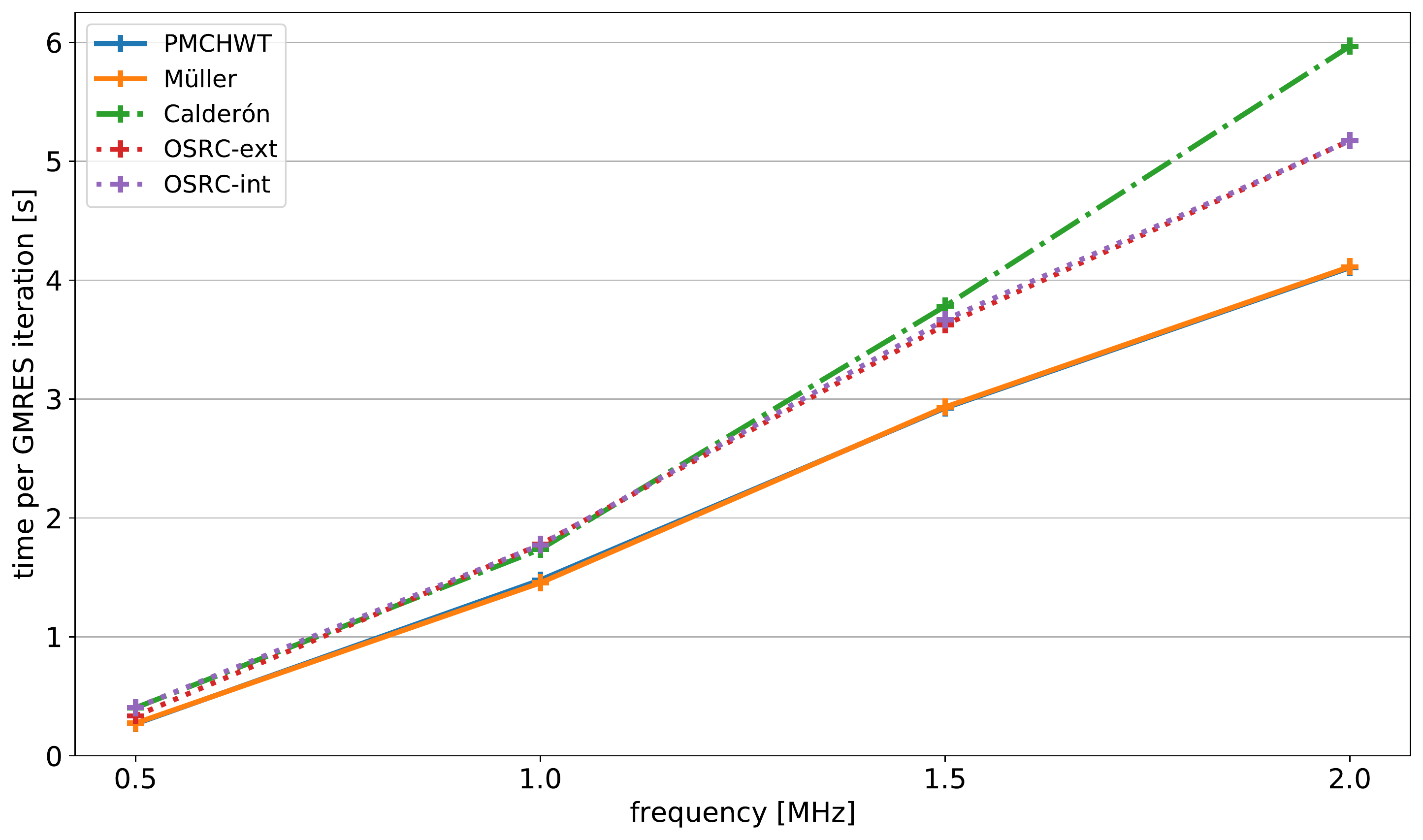}
		\caption{The computation time (wall-clock time) per GMRES iteration.}
	\end{subfigure}
	\begin{subfigure}[b]{\columnwidth}
		\centering
		\includegraphics[height=.20\paperheight]{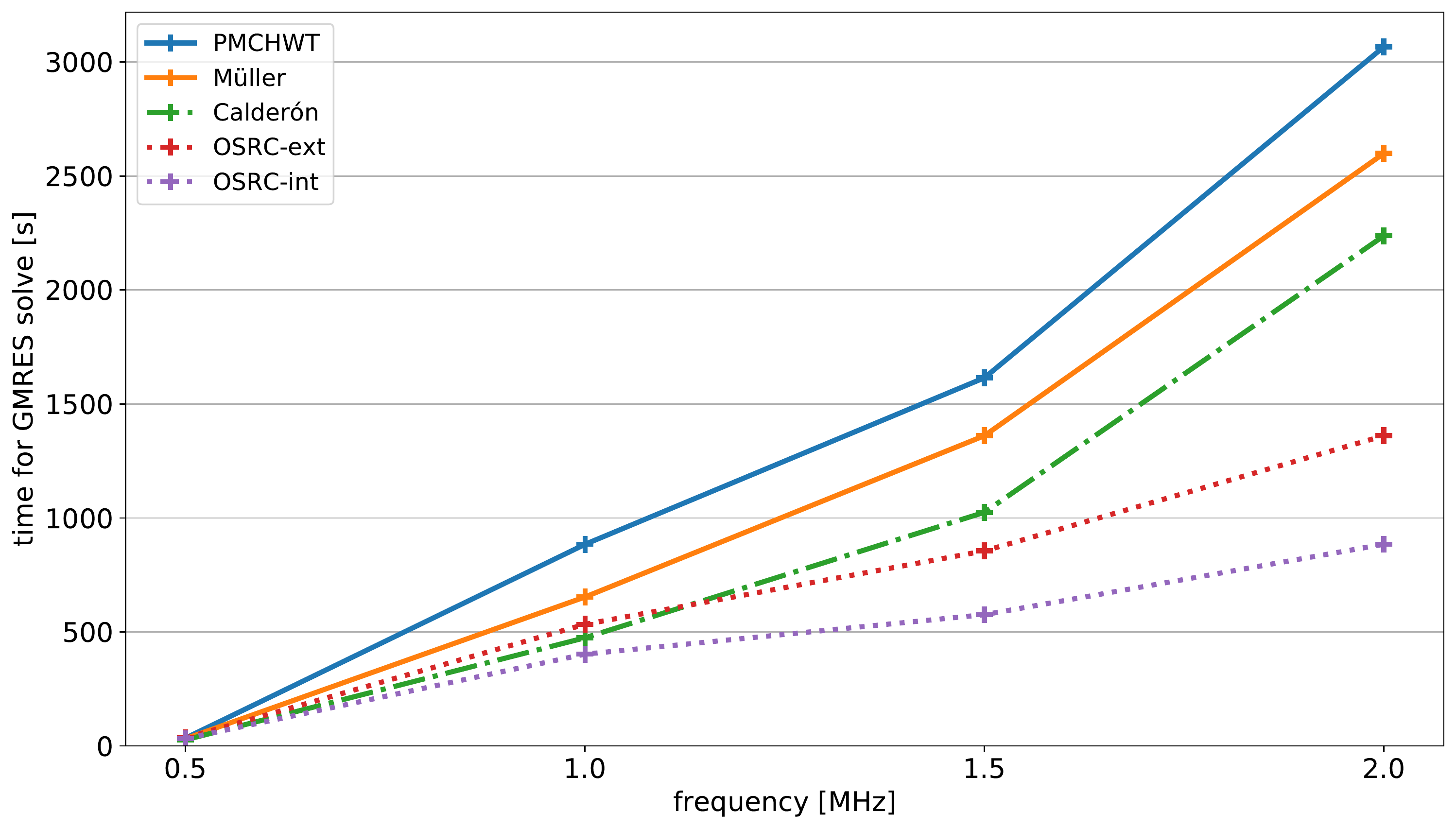}
		\caption{The total time (wall-clock time) for the GMRES algorithm.}
	\end{subfigure}
	\caption{Performance characteristics of the benchmark with four spheres.}
	\label{fig:multiple}
\end{figure}

The convergence behaviour presented in Figure~\ref{fig:multiple} indicates that the OSRC preconditioner with interior wavenumber outperforms the version with exterior parameters. That is, the wavenumber of the OSRC is based on the parameters of bone and fat at the corresponding interfaces. The Calderón preconditioner requires around half the number of iterations of the PMCHWT formulation. Since a diagonal version of the Calderón preconditioner for multiple scattering is used, the preconditioner consists of self-interaction only and the computational overhead in each GMRES iteration is only half the time for the PMCHWT formulation. The OSRC preconditioners also involve additional computation time in each GMRES iteration but because of their sparsity, the total time to solve the system is very low. The calculation time also confirms the frequency-robustness of OSRC preconditioning. That is, while Calderón preconditioning speeds up the PMCHWT formulation with a factor of 1.3 at 500~kHz and 1.4 at 2~MHz, the interior OSRC preconditioning has a speedup factor of 1.1 at 500~kHz and 3.5 at 2~MHz. Hence, the gain in efficiency improvement with OSRC preconditioning increases with frequency.

\FloatBarrier
\subsection{Large-scale acoustic scattering}

The OSRC preconditioned PMCHWT formulation is designed for high frequency acoustic scattering by penetrable objects and the computational benchmarks confirm its superior efficiency over standard boundary integral formulations when high frequencies are considered. Let us consider a large-scale problem of biomedical engineering interest to showcase the performance of the OSRC preconditioner, based on a scenario involving the treatment of an osteoid osteoma in leg bone.

An osteoid osteoma is a benign bone tumour that usually develops in the long bones of the body, such as the femur (thighbone) and tibia (shinbone). Although osteoid osteomas do not spread throughout the body, they can cause pain and discomfort. Osteoid osteomas can affect people of all ages, but they occur more frequently in children and young adults. Magnetic Resonance (MR) guided focused ultrasound has recently become one of the gold standard treatments for managing osteoid osteomas~\cite{arrigoni2020evolution}. The region on the bone to be ablated is destroyed using a focused ultrasound field, enabling the osteoid osteoma to be treated without invasive surgery. The high acoustic contrast between soft tissue and bone can lead to unwanted reflections and scattering of the focused ultrasound at the interface between soft tissue and bone. The clinical challenges ensuing from this may be addressed as part of a treatment planning stage using numerical modelling techniques such as BEM~\cite{haqshenas2021fast}.

The acoustic source modelled consists of a transducer based on the one found in the Philips Sonalleve MR-guided focused ultrasound treatment platform (Royal Philips, Amsterdam, the Netherlands). The source is in effect a 256-element array transducer, with a 140~mm geometric focal length and 140~mm diameter and can be operated at a frequency of 1.2~MHz~\cite{miloro2016feasibility}. The elements of the array are positioned along the surface of a spherical bowl in a pseudo-random fashion.

For the purpose of the BEM simulation, the transducer incident field was modelled as an array of plane circular pistons, rigidly vibrating with uniform phase and amplitude. The piston normal velocities were determined so that the array would radiate an acoustic power of 30~W in water. This resulted in a normal velocity of 0.678~m/s for each piston.

An STL file describing the femur, tibia, patella and fibula of a female right leg was obtained from GrabCAD~\cite{grabcad}. The transducer was positioned so that its focus coincided with the outer surface of the upper right side of the tibia. Given the highly focused nature of the ultrasound beam (1~mm $\times$ 1~mm $\times$ 7~mm)~\cite{philipsMR}, it is unnecessary to produce a mesh of the whole of the right leg for the purpose of the BEM simulation. Instead, a mesh was generated for only the upper section of the tibia. After editing the original STL file in Meshmixer~\cite{meshmixer}, a mesh of triangular elements was generated using Gmsh, with the following mesh statistics: 145\,338~nodes, 0.132~mm minimum element edge length, 0.501~mm maximum element edge length, and 0.295~mm average element edge length. For the purpose of this large-scale problem, the bone and transducer were immersed in water. The acoustic properties of the acoustic media are displayed in Table~\ref{table:parameters:physical}. The mesh density of the tibia section corresponds to just over four elements per wavelength of the exterior medium at a frequency of 1.2~MHz. The system matrix was compressed in a hierarchical matrix format with a tolerance of $10^{-5}$. The OSRC preconditioner uses $N_\text{Padé} = 8$, $R_\epsilon = 10^{-5}$ and the wavenumber corresponding to the interior bone material.

\begin{figure}[!ht]
	\centering
	\includegraphics[width=\columnwidth]{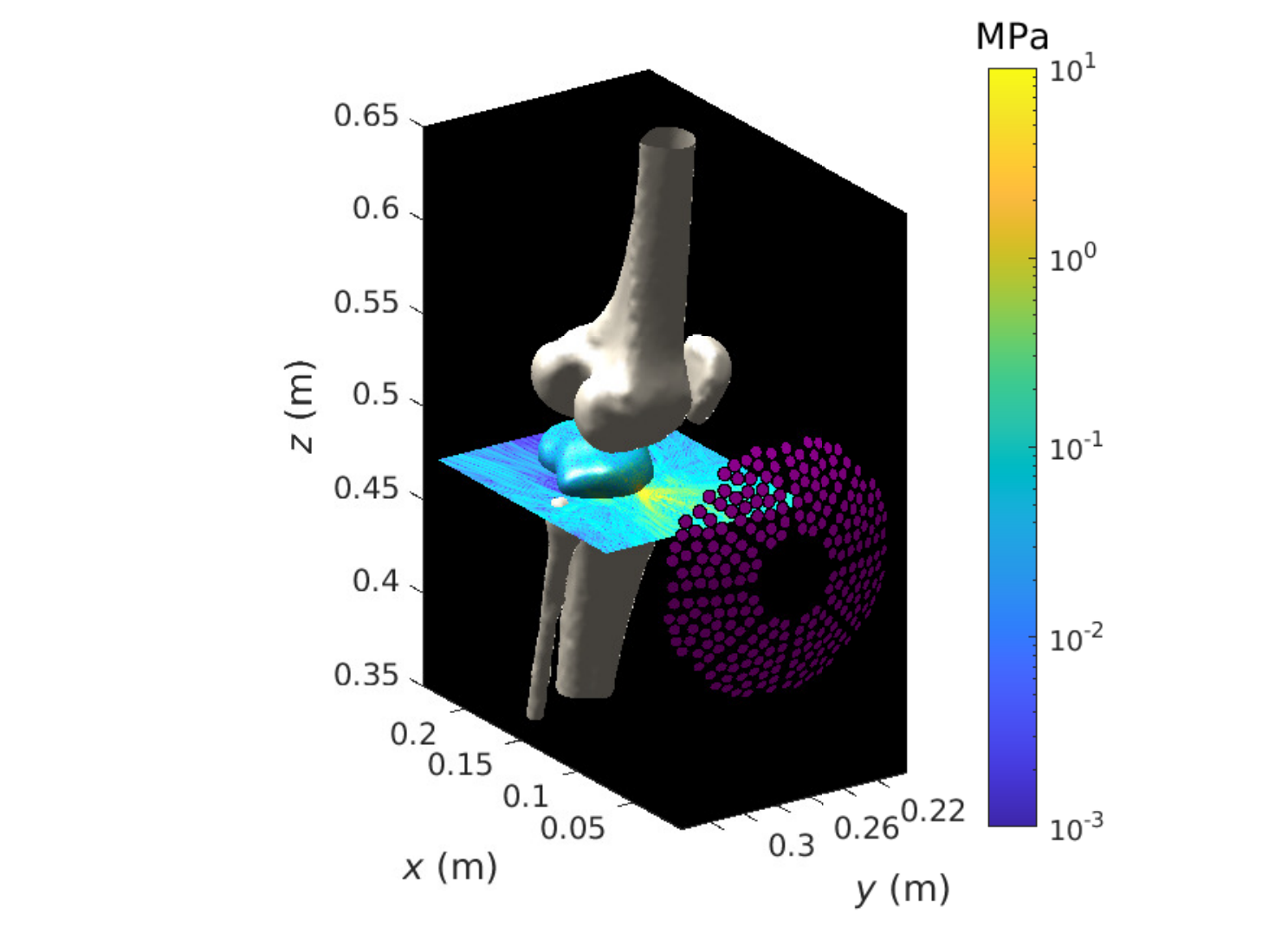}
	\caption{Acoustic pressure in focal plane and on surface of computational domain. The tibia, patella, femur and fibula are shown in bone colour. The elements of the ultrasonic transducer are shown in magenta.}
	\label{fig:hifu:setup}
\end{figure}

The system of equations was solved with GMRES in 958 iterations with an average time of 16.4~s per iteration. Figure~\ref{fig:hifu:setup} shows the position of the array transducer relative to the female right leg. The acoustic pressure in the transducer focal plane was obtained. The acoustic pressure map on the surface of the tibia section is shown in the same figure. A logarithmic colour map was used to better visualise the details of the field.

\begin{figure}[!ht]
	\centering
	\includegraphics[width=\columnwidth]{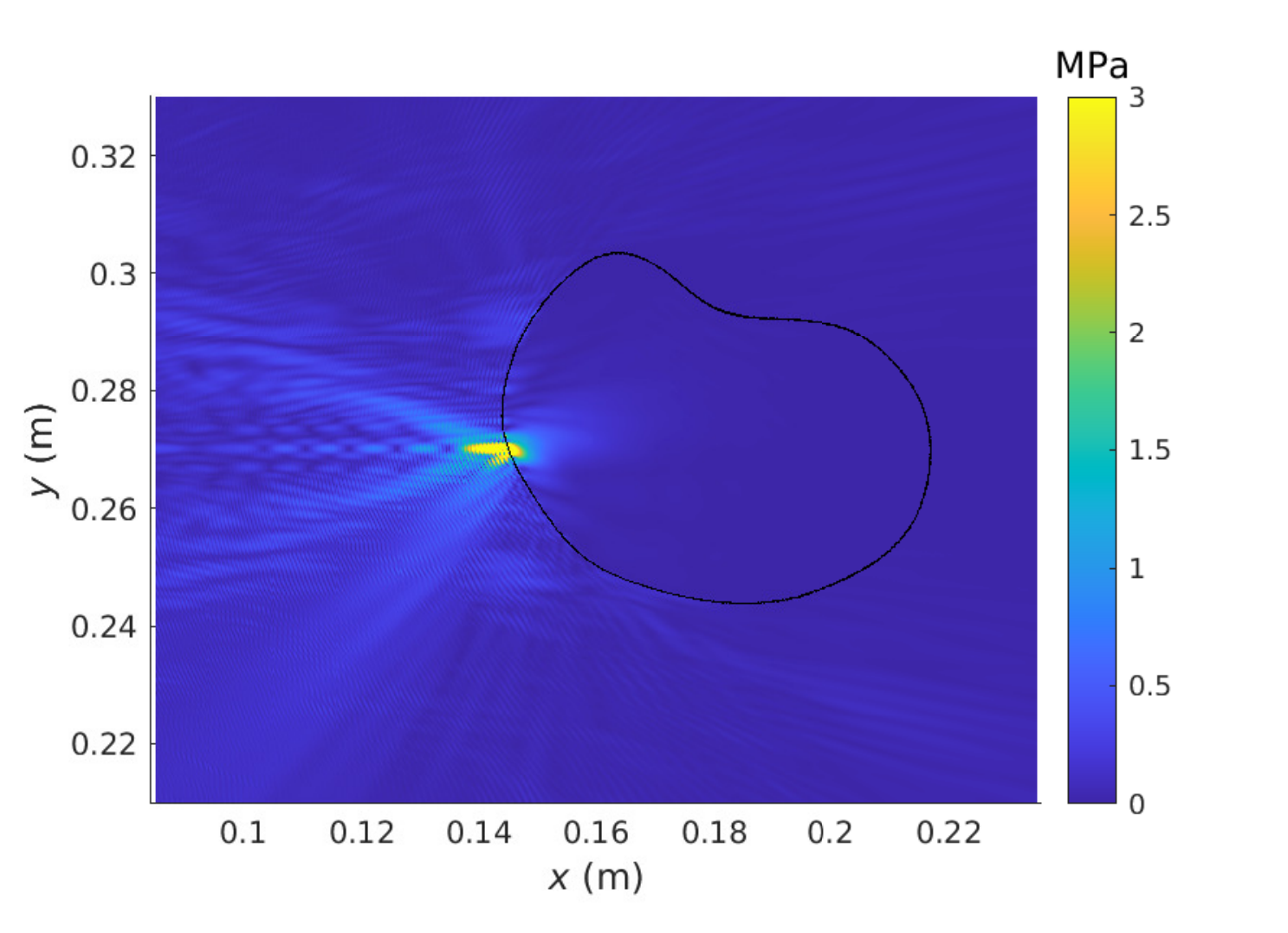}
	\caption{Acoustic pressure in the focal plane. The outer boundary of the tibia is highlighted with the black contour line.}
	\label{fig:hifu:field}
\end{figure}

The total acoustic pressure in the focal plane is displayed in Figure~\ref{fig:hifu:field}. The colour map is now linear. The focus of the transducer is positioned at $x = 140$~mm and $y = 270$~mm. Scattering by the bone is clearly visible, as is interference between the incident and reflected fields.

\FloatBarrier
\section{Conclusions}

Preconditioning is necessary to reduce the computational footprint of acoustic scattering at large-scale geometries with the BEM. This work introduced a preconditioner based on OSRC operators that effectively improves the conditioning of the PMCHWT formulation and is robust with frequency. Computational benchmarks confirm the efficiency of the frequency-robust OSRC preconditioning at canonical test cases and multiple scattering at penetrable objects. The OSRC preconditioned PMCHWT formulation can accurately simulate osteoid osteoma at a realistic model of a human knee joint at operating frequencies of ultrasound therapy.

\section*{Acknowledgment}

This work was financially supported by CONICYT [FONDECYT 11160462], the EPSRC [EP/P012434/1], and the Vicerrectoría de Investigación de la Pontifica Universidad Católica de Chile.

\bibliographystyle{elsarticle-num} 
\bibliography{refs}

\end{document}